\documentclass[letterpaper11pt]{article}

\usepackage[utf8]{inputenc}
\usepackage[T1]{fontenc}
\usepackage[breaklinks]{hyperref}
\usepackage[margin=1in]{geometry}
\usepackage{amssymb}
\usepackage{amsmath}
\usepackage{dsfont}
\usepackage{mathtools}
\usepackage{authblk}
\usepackage{xcolor}
\usepackage{dsfont}
\usepackage{subcaption}
\usepackage{wrapfig}
\usepackage{amsthm}
\usepackage{xurl}

\usepackage{microtype}
\usepackage{cite}

\hypersetup{breaklinks=true}

\newtheorem{definition}{Definition}
\newtheorem{theorem}{Theorem}
\newtheorem{proposition}{Proposition}

\newtheorem{lemma}[theorem]{Lemma}
\newtheorem*{remark}{Remark}

\newcommand{\E}{\mathbb{E}}
\newcommand{\1}{\mathds{1}}
\newcommand{\var}{\text{Var}}
\newcommand{\vard}{\text{Var}_{\Delta}}
\newcommand{\cov}{\text{Cov}}
\newcommand{\corr}{\text{Corr}}

\title{\LARGE \textbf{Control variates for variance-reduced ratio of means estimators}}
\author[1,2]{Louison Bocquet-Nouaille \thanks{Corresponding author \\ Emails: louison.bocquet\_nouaille@onera.fr, jerome.morio@onera.fr, benjamin.bobbia@isae-supaero.fr}}
\author[1,2]{Jérôme Morio}
\author[2]{Benjamin Bobbia}
\affil[1]{ONERA/DTIS, Université de Toulouse, F-31055 Toulouse}
\affil[2]{Fédération ENAC ISAE-SUPAERO ONERA, Université de Toulouse, 31000 Toulouse}
\date{\vspace{-10pt} \small{Preprint - July 7, 2026}}

\begin{document}

\title{Control variates for variance-reduced ratio of means estimators} 

\maketitle

\vspace{10pt}
The control variates method is a classical variance reduction technique for Monte Carlo estimators that exploits correlated auxiliary variables without introducing bias. In many applications, the quantity of interest can be expressed as a ratio of expectations. We propose a variance-reduced estimator for such ratios, which applies control variates to both the numerator and the denominator. The control variates coefficients are optimized jointly to minimize the approximated variance of the resulting estimator. This approach guarantees variance reduction and naturally extends to approximate control variates. Simulation studies show significant variance reduction, particularly when correlations between variables and control variates are strong. The practical value of the method is illustrated on multi-fidelity applications: estimating a proportion in an aircraft design use case and a conditional value-at-risk in an electromagnetic dataset.
\vspace{20pt}

\section{Introduction}
The variance of the Monte Carlo empirical mean estimator is inversely proportional to the number of available samples, which is often limited in real-world applications. As a result, every observation has a large influence on the estimate, leading to high variance and reduced reliability.

An overview of common variance reduction techniques can be found in \cite{asmussen_variance-reduction_2007}. In the setting considered here, where only a limited number of samples is available and further sampling is not possible, most methods are not directly applicable. The control variates method \cite[Sec.~2]{asmussen_variance-reduction_2007} is well-suited, as it exploits i.i.d. samples from an auxiliary variable correlated with the target variable, to reduce the variance of the Monte Carlo estimator without introducing bias.  It can be seen as a form of transfer learning \cite{zhu_recent_2025}, where information from a source variable is used to improve estimation on the target variable.

The control variates method was originally framed as the "regression estimator" \cite[Section 2.6]{hammersley_monte_1965} as it can be interpreted as a linear regression of the target variable on the control variate, an auxiliary variable. The optimal regression coefficient is the slope that best explains the variability of the target variable through its correlation with the control variate. The following selection of works highlights the method's adaptation to various mathematical and practical settings. In terms of functional form, \cite{lewis_smarter_1987} introduced both linear and nonlinear alternatives to the classical additive form of the control variates estimator. From a theoretical perspective, the method has been formalized within a comprehensive Hilbert-space framework \cite{glynn_new_2002}, establishing asymptotic optimality results. In the context of Monte Carlo integration, the method has been augmented by modern approximation techniques, ranging from classical regression-based approaches \cite{salaun_regression-based_2022} to flexible frameworks like neural networks \cite{muller_neural_2020} and polynomial chaos expansions \cite{fox_polynomial_2021}. Furthermore, the control variates method has been adapted for quasi-Monte Carlo estimators \cite{hickernell_control_2005}, constrained variable estimation \cite{maskell_control_2022}, and Markov chain Monte Carlo (MCMC) methods, to be applied on dependent samples, as in \cite{dellaportas_control_2012}.
The methodology has also been tailored to specific estimation targets, with works by \cite{hesterberg_control_1998} and \cite{davidson_regression-based_1992} proposing specific control variate functions for probabilities and quantiles. Recent advances have moved toward a more general, functional perspective. \cite{oates_control_2017, south_semi-exact_2022} extend the control variates methodology into the domain of Control Functionals (CFs), which leverage Stein’s method and log-posterior derivatives to automatically construct optimal functional operators for variance reduction. This framework eliminates the need for predefined auxiliary variables, and has recently been extended to vector-valued settings \cite{sun_vector-valued_2023}.

The approximate control variates method extends the standard "exact" control variates approach to cases where the mean of the auxiliary variable is unknown, estimating it from additional samples. It is also referred to as Multi-Fidelity Monte Carlo (MFMC) \cite{cannamela_controlled_2008, gorodetsky_generalized_2020, kim_parametric_2024}: the output of a computationally cheaper lower-fidelity model is used as a control variate to improve the estimation on the output of a more accurate high-fidelity model. Since both models approximate the same underlying quantity, their outputs are often highly correlated, which enables efficient variance reduction through control variates. \\

In this work, we focus on ratio of means estimation for real random variables. Among numerous applications, the ratio of means estimator is a quantity of interest with various structural roles across a wide range of disciplines. It is used to estimate proportions, such as in aircraft design \cite{meheut_strut_braced_2024} and epidemiology via incidence and mortality rates \cite{ferlay_estimates_2007}. It also naturally arises when estimating conditional expectations, as seen in reliability analysis via the Conditional Value-at-Risk (CVAR) \cite{hong_monte_2014} and self-normalized importance sampling \cite{bugallo_adaptive_2017}. Finally, the estimator is well-suited for financial ratio analysis \cite{lin_financial_2011}. These examples are illustrative rather than exhaustive, highlighting the ratio of means estimator's ubiquity across diverse fields.

The standard approach to estimate ratios of means consists of estimating the numerator and denominator separately by Monte Carlo, and then taking their ratio. This estimator has been analyzed in \cite[Chapter 6]{cochran_sampling_1977}. The control variates method can be applied to both components of the ratio of means estimator to reduce variance. To the best of our knowledge, \cite{gordon_efficient_1982} is the only work specifically extending the method to ratio of means estimation, by setting the numerator coefficient to its classical value, the one that minimizes the variance of a single mean, and optimizing the denominator coefficient to minimize the approximated ratio variance.

The classical coefficients for mean estimation, or those of \cite{gordon_efficient_1982}, can in fact increase the variance of the ratio of means estimator with certain correlation structures, as they are not both derived to minimize the variance of the ratio. Because these methods optimize coefficients independently or sequentially, they operate within a limited parameter subspace that may exclude the global minimum. This deviation from the optimal control variates coefficients can result in a significant loss of efficiency. This paper introduces new optimal coefficients that guarantee variance reduction as they are jointly optimized to minimize the estimator's variance, approximated by the Delta method. The proposed closed-form optimal coefficients extend the classical variance reduction guarantees of control variates from simple means to ratio of means estimators. \\

The paper is organized as follows. In Section \ref{sec:CV}, the control variates method is presented for mean and ratio of means estimation, covering both exact and approximate forms for each. Section \ref{sec:contrib} introduces the new optimal coefficients for variance-reduced ratio of means estimators using control variates, along with an elliptical stability condition for variance reduction. In Section \ref{sec:simu}, simulation results are presented, illustrating extreme correlation structures and their impact on variance reduction. Finally, in Section \ref{sec:appli}, we demonstrate the proposed estimators on two multi-fidelity datasets: to estimate a proportion in an aircraft design dataset, and a Conditional Value-at-Risk on an electromagnetic dataset. For readability, all proofs are collected in the supplementary material, in Appendix \ref{sec:proofs}.

\section{Variance reduction with control variates} \label{sec:CV}

\subsection{Control variates applied to a mean estimator}
\subsubsection{Exact control variates}
The control variates method aims to reduce the variance of the Monte Carlo estimator defined as follows. 
\begin{definition}[Monte Carlo estimator] \label{def:MC}
Let $(A_i)_{i=1 \dots n}$ be i.i.d. samples of a random variable $A \in L^2$. The Monte Carlo estimator $\overline{A_n}$ of $\E[A]$ is defined as \begin{align}\overline{A_n}=\frac{1}{n}\sum_{i=1}^{n}A_i \label{eq:MC}.\end{align}
\end{definition}
The variance of the Monte Carlo estimator is proportional to the number of samples $n$:\begin{align} \var(\overline{A_n}) = \frac{1}{n} \var(A). \end{align} To reduce the variance without increasing $n$, the idea is to consider a control variate $B$ correlated with the target random variable $A$, and with known expectation $\E[B]$. This approach is a form of statistical transfer learning: it leverages knowledge from a related variable $B$ to compensate for the lack of information about the target variable $A$. The control variates estimator of $\E[A]$ is defined as follows. 
\begin{definition}[Control variates estimator \cite{asmussen_variance-reduction_2007}] \label{def:CV}
Let $(A_i, B_i)_{i=1 \dots n}$ be i.i.d. samples from the joint distribution of the random variables $A, B \in L^2$. Assume $\E[B]$ is known. The control variates estimator $\widehat{A}_{CV}$ of $\E[A]$ is defined as \begin{align}\widehat{A}_{CV}(\alpha_c) = \overline{A_n} + \alpha_c (\E[B] - \overline{B_n}) ,\end{align} where the coefficient $\alpha_c$ is set to minimize the variance of the estimator: \begin{align}\alpha_c := \underset{\alpha \in \mathbb{R}}{\text{argmin}} \; \var\left(\overline{A_n} + \alpha (\E[B] - \overline{B_n})\right) = \frac{\cov(A, B)}{\var(B)}. \label{eq:CV}\end{align}
\end{definition}

With the optimal coefficient $\alpha_c$, the variance of the control variates estimator is \begin{align} \var(\widehat{A}_{CV}(\alpha_c)) = \left(1 - \corr(A,B)^2\right) \var(\overline{A_n}) ,\end{align} where $\corr(A,B)$ is the linear correlation coefficient between $A$ and $B$. This expression shows that the method never increases variance. The extent of variance reduction depends on the strength of the correlation between $A$ and $B$: the higher the absolute value of the correlation, the greater the reduction. When the correlation is zero, the control variates estimator coincides with the Monte Carlo one.

The coefficient $\alpha_c$ needs to be estimated, which can introduce a small bias if the estimation is made from the same dataset used to compute $\widehat{A}_{CV}(\alpha_c)$. This plug-in bias \cite{kennedy2024semiparametric} is generally negligible in practice, as discussed in \cite{owen_variance_2013}. Note that if $\alpha_c$ is so poorly estimated that $\widehat{\alpha}_c \notin (\min(0, 2\alpha_c), \max(0, 2\alpha_c))$, the variance will increase. The proof of this condition is presented in Appendix \ref{sec:cond_est_coef_cv}.

\begin{remark}
    In practice, the control variate $B$ is often designed to "mirror" the target variable $A$ by evaluating the same function defining $A$ on an auxiliary dataset \cite{hesterberg_control_1998}. While intuitive, this structural choice is not necessarily the best approach as the variance reduction achieved is strictly limited by how effectively the construction of $B$ maximizes the correlation with $A$ for the given sample sets.
\end{remark}
Control Functionals \cite{oates_control_2017, south_semi-exact_2022} offer an automated alternative to manual selection of the control variate. However, we restrict our scope here to the classical regression framework in which the control variate is defined beforehand.

While the method naturally extends to multiple control variates using linear regression to estimate the various coefficients \cite{leluc_control_2021}, only a single control variate is considered here for clarity.

\subsubsection{Approximate control variates}
In practice, the assumption that the control variate mean $\E[B]$ is known is often unrealistic. The approximate control variates (ACV) method overcomes this by estimating the mean from additional data. It is also referred to as Multi-Fidelity Monte Carlo (MFMC) \cite{cannamela_controlled_2008, gorodetsky_generalized_2020, kim_parametric_2024}. In this framework, the output of a computationally cheaper, lower-fidelity model is used as a control variate to improve the estimation of a more accurate high-fidelity model.
\begin{definition}[Approximate control variates \cite{gorodetsky_generalized_2020}]
Let $(A_i, B_i)_{i=1 \dots n}$ be i.i.d. samples from the joint distribution of the random variables $A, B \in L^2$, and let $(B_i)_{i=n+1 \dots n+m}$ be additional i.i.d. samples of $B$. The approximate control variates estimator of $\E[A]$ is defined as \begin{align}\widehat{A}_{ACV}(\alpha_c) = \overline{A_n} + \alpha_c (\overline{B_{n+m}} - \overline{B_n}).\end{align}
\end{definition}
The optimal coefficient $\alpha_c$ coincides with the one derived for the exact case, given by \eqref{eq:CV}.

The variance of the approximate control variates estimator is \begin{align} \var(\widehat{A}_{ACV}(\alpha_c)) = \left(1 - \frac{m}{n+m} \corr(A,B)^2\right) \var(\overline{A_n}). \label{eq:approx_red}\end{align} This expression highlights that increasing the number of additional samples from the control variate $B$ leads to greater variance reduction. For a fixed $n$, as $m \to \infty$, the estimator variance converges to the variance of the exact control variates estimator, as $\overline{B_{n+m}} \overunderset{a.s.}{m \to \infty}{\to} \E[B]$. \\

\subsection{Control variates applied to ratio of means estimators}
\subsubsection{Exact control variates}
This section presents a variance-reduced estimator based on control variates, for the ratio of means $R$ defined as \begin{align}R = \frac{\E[A]}{\E[C]},\end{align} where $A$ and $C$ are arbitrary random variables. Assume $\E[C] \neq 0$. The control variates method can be applied to reduce the variance of the Monte Carlo ratio of means estimator \begin{align}\widehat{R}_{\frac{MC}{MC}} = \frac{\overline{A_n}}{\overline{C_n}},\end{align} referred to as the MC/MC estimator, where $\overline{A_n}$ and $\overline{C_n}$ are the Monte Carlo estimators of $\E[A]$ and $\E[C]$, as in \eqref{eq:MC}. 

The MC/MC estimator is asymptotically unbiased as $n$ goes to infinity. For smaller sample sizes, the nonlinearity of the estimator introduces a bias, as the expectation of a ratio is generally not equal to the ratio of expectations. This bias is particularly relevant when the denominator has high variance or is near zero \cite[Section 6.5]{cochran_sampling_1977}. \\

To reduce the variance of the MC/MC estimator, the control variates method can be applied to both the numerator and the denominator, with two control variates $B$ and $D$ with known means, used to "control" $A$ and $C$ respectively. We refer to the resulting estimator as the CV/CV estimator. 
\begin{definition}[CV/CV estimator]
Let $(A_i, B_i, C_i, D_i)_{i=1 \dots n}$ be i.i.d. samples from the joint distribution of the random variables $A, B, C, D \in L^2$, where $B$ and $D$ are control variates. Assume $\E[B]$ and $\E[D]$ are known. The CV/CV estimator of $R$ is defined as \begin{align}\widehat{R}_{\frac{CV}{CV}}\left(\alpha,\beta\right) = \frac{\widehat{A}_{CV}(\alpha)}{\widehat{C}_{CV}(\beta)} = \frac{\overline{A_n} + \alpha (\E[B] - \overline{B_n})}{\overline{C_n} + \beta (\E[D] - \overline{D_n})}. \end{align}
\end{definition}

The asymptotic properties of the MC/MC estimator extend to the control variates estimator $\widehat{R}_{\frac{CV}{CV}}$ defined previously, notably the absence of asymptotic bias. As illustrated in Sections \ref{sec:simu} and \ref{sec:appli}, the impact of the bias in finite samples is negligible in practice. Given that both bias and variance scale with $1/n$, the estimator remains highly robust even for moderate sample sizes.\\

Let $\vard$ denote the linear approximation of the variance via the Delta method when $n \to \infty$, which enables the closed-form derivation of control variates coefficients. The explicit expansion \eqref{eq:var_mc_mc} is provided in Appendix \ref{sec:proof_var_mc_mc}. Its residual error is of order $\mathcal{O}_{\infty}(n^{-3/2})$.

\begin{proposition}[Coefficients for the CV/CV estimator]
Different choices of coefficients exist:
    \begin{itemize}
        \item the classical coefficients for mean estimators, introduced in Definition \ref{def:CV}: \begin{alignat}{2}\alpha_c := \underset{\alpha \in \mathbb{R}}{\text{argmin}} \; \var\left(\overline{A_n} + \alpha (\E[B] - \overline{B_n})\right) = \frac{\cov(A, B)}{\var(B)}; \label{eq:classic_alpha} \\
        \beta_c := \underset{\beta \in \mathbb{R}}{\text{argmin}} \; \var\left(\overline{C_n} + \beta (\E[D] - \overline{D_n})\right) = \frac{\cov(C, D)}{\var(D)}; \label{eq:classic_beta}\end{alignat} 
        \item the coefficients introduced in \cite{gordon_efficient_1982}: \begin{alignat}{2}
        \alpha_g &:= \underset{\alpha \in \mathbb{R}}{\text{argmin}} \; \var\left(\overline{A_n} + \alpha (\E[B] - \overline{B_n})\right) = \frac{\cov(A, B)}{\var(B)}; \label{eq:gordon_alpha} \\
        \beta_g  &:= \underset{\beta \in \mathbb{R}}{\text{argmin}} \; \vard\left(\frac{\overline{A_n} + \alpha_g (\E[B] - \overline{B_n})}{\overline{C_n} + \beta (\E[D] - \overline{D_n})}\right) = \frac{R\cov(C, D) - \cov(A,D) + \alpha_g \cov(B,D)}{R\var(D)}. \label{eq:gordon_beta}
        \end{alignat}
    \end{itemize}
\end{proposition}

However, these choices are not always optimal: variance reduction with respect to the MC/MC estimator is not guaranteed in all cases, and the variance can even be increased as illustrated with simulations in Section \ref{sec:simu}. The variance difference between the CV/CV and MC/MC estimators can be explicitly expressed as a function of the covariances and variances of $A$, $B$, $C$, and $D$ (see Appendix \ref{sec:diff_var}). Crucially, Theorem \ref{prop:cond_est_coef} shows that the CV/CV estimator outperforms the MC/MC one if and only if the estimated coefficients fall within a specific ellipse, centered at the optimal coefficients defined in Proposition \ref{prop_cv_cv}.

\subsubsection{Approximate control variates}
Assuming that the control variates' means $\E[B]$ and $\E[D]$ are known is often an unrealistic assumption, which can be overcome by the approximate control variates method, which estimates the means from additional data.
\begin{definition}[ACV/ACV estimator]
Let $(A_i, B_i, C_i, D_i)_{i=1 \dots n}$ be i.i.d. samples from the joint distribution of the random variables $A, B, C, D \in L^2$, and let $(B_i,D_i)_{i=n+1 \dots n+m}$ be additional i.i.d. samples of $B$ and $D$. The ACV/ACV estimator of $R$ is defined as \begin{align}
    \widehat{R}_{\frac{ACV}{ACV}}\left(\alpha,\beta\right) = \frac{\widehat{A}_{ACV}(\alpha)}{\widehat{C}_{ACV}(\beta)} = \frac{\overline{A_n} + \alpha (\overline{B_{n+m}} - \overline{B_n})}{\overline{C_n} + \beta (\overline{D_{n+m}} - \overline{D_n})}.
\end{align}
\end{definition}

\section{Optimal coefficients and associated properties} \label{sec:contrib}
\subsection{Optimal coefficients}
In this section, we propose a new set of optimal coefficients that guarantee asymptotic variance reduction with the CV/CV and ACV/ACV estimators. While existing methods optimize the control variates for the numerator and denominator independently or sequentially, the proposed approach relies on a joint optimization of the coefficients, minimizing the approximated variance of the ratio \eqref{eq:var_cv_cv} with residual error of order $\mathcal{O}_{\infty}(n^{-3/2})$.
\begin{proposition}[Optimal coefficients for the CV/CV estimator] \label{prop_cv_cv} 
Assume $|\corr(B,D)|<1$. The optimal coefficients for the CV/CV estimator are defined as:
    \begin{align*}
    (\alpha_o, \beta_o) := \underset{(\alpha,\beta) \in \mathbb{R}^2}{\text{argmin}} \; \vard\left(\frac{\overline{A_n} + \alpha (\E[B] - \overline{B_n})}{\overline{C_n} + \beta (\E[D] - \overline{D_n})}\right).
    \end{align*}
The closed-form expressions for $\alpha_o$ and $\beta_o$ are given by:
    \begin{alignat}{2}
    \alpha_o = \frac{\var(D)\cov(A, B) - R \var(D)\cov(B, C)  + R \cov(B, D)\cov(C, D) -  \cov(B, D)\cov(A, D)}{\var(B)\var(D) - \cov(B,D)^2} ; \label{eq:alpha_opt} \\
    \beta_o = \frac{ \cov(B, D)\cov(A, B) -  R \cov(B, D) \cov(B, C) + R \var(B)\cov(C, D) - \var(B)\cov(A, D)}{R \left(\var(B)\var(D) - \cov(B,D)^2 \right)} , \label{eq:beta_opt}  \end{alignat}
where $R$ is the ratio to estimate.
\end{proposition}

The full derivation of the optimal coefficients $\alpha_o$ and $\beta_o$ is presented in Appendix \ref{sec:cv_cv_opt_coef}. It is based on an approximated expression of the CV/CV estimator’s variance, which is shown to be a convex function of $(\alpha, \beta)$. The closed-form formulas of the optimal coefficients are then obtained by setting the gradient of this variance to zero and solving for $(\alpha, \beta)$.

The optimal coefficients of the ACV/ACV estimator are the same as the ones of the CV/CV estimator, as illustrated in Appendix \ref{sec:approx}.

The constraint $|\corr(B, D)| < 1$ is essential to ensure that the approximated variance of the ratio is a strictly convex function of $(\alpha, \beta)$. If this condition is violated, the optimal coefficients are no longer uniquely defined. From a practical standpoint, if $|\corr(B, D)| \to 1$, the empirical estimators of $\alpha_o$ and $\beta_o$ exhibit high variance. The resulting estimation error can lead to a variance increase in the final CV/CV estimator, a phenomenon further detailed in Theorem \ref{prop:cond_est_coef}.

If the control variates $B$ and $D$ are chosen such that $|\corr(B,D)| = 1$, alternative variance-minimizing coefficients $\Tilde{\alpha_o}$ and $\Tilde{\beta_o}$ are given in Proposition \ref{prop_B=D}, in Appendix \ref{sec:prop_B=D}. However, they no longer guarantee variance reduction, which is instead contingent upon the inequality given in \eqref{cond}.

\subsection{Asymptotic variance reduction guarantee}
\begin{proposition}
    With the optimal coefficients given in \eqref{eq:alpha_opt} and \eqref{eq:beta_opt}, the approximated variance reduction obtained with the CV/CV estimator is: \begin{align}
    &\vard\left(\widehat{R}_{\frac{CV}{CV}}\left(\alpha_o,\beta_o\right)\right) - \vard\left(\widehat{R}_{\frac{MC}{MC}}\right) \nonumber\\ &= - \frac{1}{n\E[C]^2} \frac{Var\bigg(\big(R\cov(B,C)-\cov(A,B)\big)D - \big(R\cov(C,D)-\cov(A,D)\big)B\bigg)}{\var(B)\var(D)-\cov(B,D)^2} + \mathcal{O}_{\infty}(n^{-3/2}) \leq 0. \label{eq:var_red}
    \end{align}
    The variance reduction is guaranteed, extending the property of standard control variates for simple means to ratio of means estimators.
\end{proposition}
The detailed derivation is presented in Appendix \ref{sec:diff_var}.

The variance difference formula given in \eqref{eq:var_red} indicates that many different covariance structures lead to variance reduction. Intuitively, one might maximize the correlations between the target variables and their respective control variates. It is, however, not straightforward to determine which covariances and variances should be maximized or minimized to boost variance reduction, as it depends on their signs.

A similar approximated variance reduction formula holds for the ACV/ACV estimator, which scales the exact CV/CV case by a factor of $\frac{m}{n+m}$: \begin{align}
    \vard\left(\widehat{R}_{\frac{ACV}{ACV}}\left(\alpha,\beta\right)\right) - \vard\left(\widehat{R}_{\frac{MC}{MC}}\right) = \frac{m}{n+m} \left( \vard\left(\widehat{R}_{\frac{CV}{CV}}\left(\alpha,\beta\right)\right) - \vard\left(\widehat{R}_{\frac{MC}{MC}}\right) \right).
\end{align}
As the optimal coefficients remain identical, variance reduction is also guaranteed for the ACV/ACV estimator.

\subsection{Stability condition with estimated coefficients}
The following theorem establishes the necessary condition for achieving asymptotic variance reduction with estimated coefficients.
\begin{theorem} \label{prop:cond_est_coef}
Let $R = \E[A]/\E[C]$ be the ratio to be estimated. Let $x_o = (\alpha_o, R\beta_o)^T$ denote the optimal coefficient vector, with $\alpha_o$ and $\beta_o$ defined in \eqref{eq:alpha_opt} and \eqref{eq:beta_opt}, and let $\widehat{x} = (\widehat{\alpha}, R\widehat{\beta})^T$ be a vector of estimated coefficients. Define the covariance matrix:$$\Lambda = \begin{pmatrix} \var(B) & -\cov(B,D) \\ -\cov(B,D) & \var(D) \end{pmatrix}.$$ The CV/CV estimator yields a lower approximated variance than the MC/MC estimator if and only if the estimated coefficients $\widehat{\alpha}$ and $\widehat{\beta}$ fall within a specific ellipse centered at the optimum $x_o$:
\begin{align}
    \vard\left(\widehat{R}_{\frac{CV}{CV}}\left(\widehat{\alpha}, \widehat{\beta}\right)\right) < \vard\left(\widehat{R}_{\frac{MC}{MC}}\right) \; \Longleftrightarrow \; (\widehat{x}-x_o)^T \Lambda (\widehat{x}-x_o) < x_o^T \Lambda x_o \, .
\end{align}
\end{theorem}
The derivation is presented in Appendix \ref{sec:cond_est_coef_cv_cv}.\\
Crucially, Theorem \ref{prop:cond_est_coef} explains why previously proposed coefficients may inadvertently increase variance, as their theoretical values can lie outside the safety zone, even before estimation noise is introduced.\\
The geometry of this ellipse is governed by the spectral decomposition of $\Lambda$. In cases where the control variates exhibit high mutual correlation, the covariance matrix becomes ill-conditioned, resulting in an elongated stability ellipse. This implies a heightened sensitivity to estimation errors in the direction of the principal eigenvector, where even minor deviations in $\hat{x}$ can push the estimator outside the stable region, leading to a variance that exceeds the baseline MC/MC estimator.

To mitigate the risks, one can alternatively use what we refer to as the CV/MC estimator, known in the literature as the top-controlled estimator \cite{bauer1987control} (defined in Appendix \ref{sec:def_cv_mc}), which applies control variates only to the numerator. This estimator requires estimating just a single optimal coefficient, providing a more reliable option when the sample size $n$ is small, or when the control variates $B$ and $D$ are highly correlated.

\section{Simulation study} \label{sec:simu}
The reported simulation values are averaged over 10,000 independent repetitions, comparing different coefficients. 
For the CV/CV estimator, three sets of coefficients are compared:
\begin{itemize}
    \item the "classical" coefficients $(\alpha_c, \beta_c)$, given in \eqref{eq:classic_alpha} and \eqref{eq:classic_beta}, derived by optimizing the numerator and denominator independently;
    \item the coefficients proposed by \cite{gordon_efficient_1982} $(\alpha_g, \beta_g)$, given in \eqref{eq:gordon_alpha} and \eqref{eq:gordon_beta}, where $\alpha_g$ minimizes the numerator variance and $\beta_g$ minimizes the approximated variance of the ratio;
    \item the new optimal coefficients $(\alpha_o, \beta_o)$, given in \eqref{eq:alpha_opt} and \eqref{eq:beta_opt}, that jointly minimize the approximated variance of the ratio.
\end{itemize}

Performance is measured with the relative variance reduction (RVR), that quantifies the variance reduction achieved by $\widehat{R}_{\frac{CV}{CV}}\left(\alpha,\beta\right)$ compared to the baseline estimator $\widehat{R}_{\frac{MC}{MC}}$: \begin{align}\text{RVR} := \frac{\var\left(\widehat{R}_{\frac{MC}{MC}}\right) - \var\left(\widehat{R}_{\frac{CV}{CV}}\left(\alpha,\beta\right)\right)}{\var\left(\widehat{R}_{\frac{MC}{MC}}\right)}. \end{align}
A positive RVR indicates that the estimator $\widehat{R}_{\frac{CV}{CV}}\left(\alpha,\beta\right)$ achieves a lower variance than the baseline $\widehat{R}_{\frac{MC}{MC}}$, while a negative RVR implies a variance increase which may occur in practice when coefficients are either sub-optimally defined or poorly estimated from finite samples. It is estimated by measuring the variance of the estimators across the 10,000 independent repetitions.\\

To assess the validity of the variance approximation used to derive the proposed optimal coefficients $(\alpha_o, \beta_o)$, we compute the Relative Approximation Error (RAE). This metric compares the empirical variance $\widehat{\var}$, observed across 10,000 independent datasets, to the approximated variance $\vard$. The latter is computed by averaging across those same 10,000 realizations the $\mathcal{O}_{\infty}(n^{-3/2})$ Taylor expansion defined in Appendix \ref{sec:var_cv_cv}. The RAE is defined as:
$$
\text{RAE} := \frac{\left| \vard - \widehat{\var} \right|}{\widehat{\var}}.
$$
Low RAE values across our experiments confirm that the first-order expansion remains a reliable proxy for the true estimator variance in the regimes considered.\\

The random variables $A, B, C,$ and $D$ are modeled as components of a four-dimensional Gaussian vector $(A, B, C, D)^\top \sim \mathcal{N}\left( \mu, \Sigma \right)$ with \begin{align*}
\mu = \begin{pmatrix} 50 \\ 20 \\ 10 \\ 100
\end{pmatrix} \; \text{and} \;  \Sigma = \begin{pmatrix}
1 & \cov(A,B) & \cov(A,C) & \cov(A,D) \\
\cov(A,B) & 1 & \cov(B,C) & \cov(B,D) \\
\cov(A,C) & \cov(B, C) & 1 & \cov(C,D) \\
\cov(A,D) & \cov(B,D) & \cov(C,D) & 1
\end{pmatrix}.
\end{align*} With these parameters, the ratio $R=\E[A]/\E[C]$ to be estimated equals 5. All variances are normalized to one, so the covariance matrix coincides with the linear correlation matrix. For clarity, only the lower triangle of the covariance matrix is reported in each experiment.\\
The sample size $n$ is 100. \\

To assess how coefficients influence the CV/CV estimator, we compare performance across three covariance structures, defined as the set of pairwise covariances between $(A, B, C, D)$. The estimator's efficiency is directly driven by these parameters as derived in the variance formulas in Appendix \ref{sec:diff_var}, and as detailed in the condition for variance reduction with estimated coefficients in Theorem \ref{prop:cond_est_coef}. The first covariance structure corresponds to the best-case scenario for the optimal coefficients proposed in this work, highlighting the maximum potential variance reduction. The second covariance structure corresponds to the worst-case scenario for the coefficients of \cite{gordon_efficient_1982}, demonstrating the maximum variance increase possible. The last covariance structure corresponds to the best-case scenario for the coefficients of \cite{gordon_efficient_1982}, providing a fair comparison between coefficients. The different covariance configurations are obtained by maximizing or minimizing the variance reduction using the approximated formulas in Appendix \ref{sec:diff_var}. The optimization is carried out using the differential evolution algorithm \cite{storn_differential_1997}, a robust global optimizer suitable for the non-convex nature of the problem, with the objective function set to infinity whenever the tested covariance matrix is not positive definite.

We emphasize that the presented covariance matrices are illustrative examples rather than the only scenarios that produce the following results. They are not 'necessary' conditions; for example, changing which variables are uncorrelated or flipping the signs in the matrix may lead to similar results. We use these scenarios to show how the estimator performs under specific regimes, acknowledging that many different structures can trigger the same statistical behavior.

\subsection{Strong variance reduction} \label{sec:simu_MC_CV}
The presented results are the best-case scenario for the proposed optimal coefficients, observed with the following covariance matrix $$ \Sigma = 
\begin{pmatrix}
1 & & & \\
0.17 & 1 & & \\
-0.99 & -0.17 & 1 & \\
0.94 & -0.14 & -0.95 & 1
\end{pmatrix}.
$$

\begin{figure}[!h]
    \centering
    \begin{subfigure}[t]{0.49\linewidth}
        \centering
        \includegraphics[width=\linewidth]{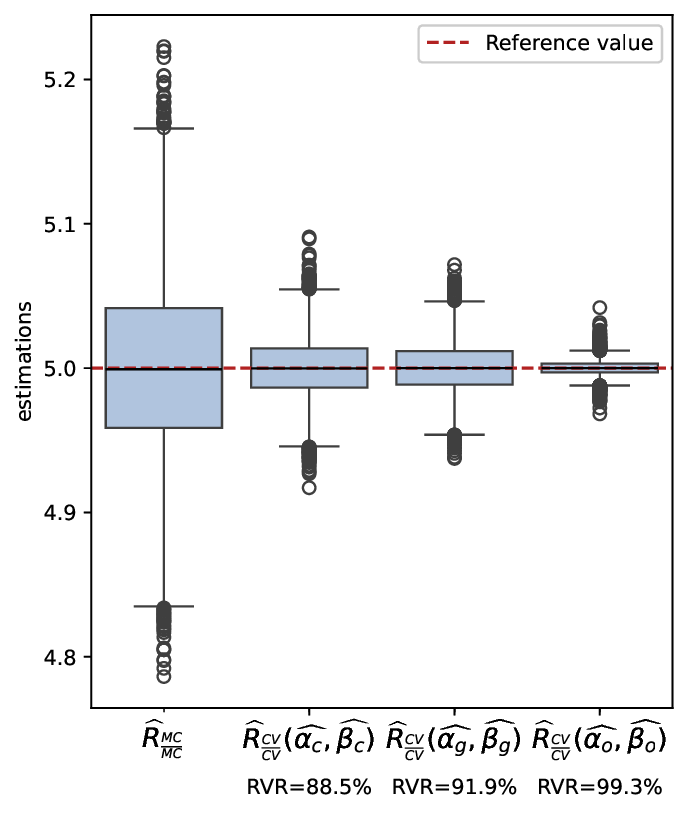}
        \subcaption{Boxplots of ratio estimations.}
        \label{fig:boxplot_all_good}
    \end{subfigure}
    \hfill
    \begin{subfigure}[t]{0.49\linewidth}
        \centering
        \includegraphics[width=\linewidth]{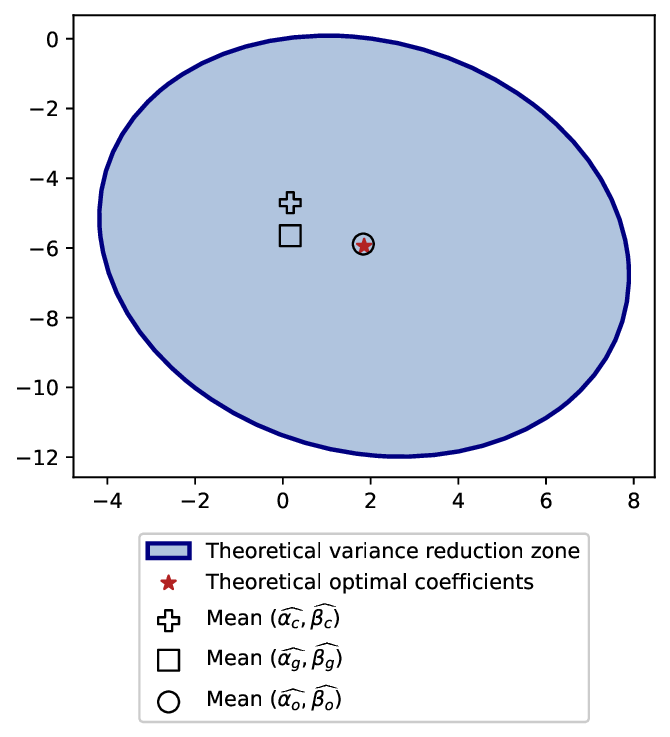}
        \subcaption{Stability ellipse and mean coefficient estimates (Theorem \ref{prop:cond_est_coef}).}
        \label{fig:ellipse_all_good}
    \end{subfigure}
    \caption{\textbf{Best-case scenario for optimal coefficients.} All coefficients achieve strong variance reduction.}
\end{figure}

In this scenario, all coefficients achieve variance reduction as evidenced in Figure \ref{fig:boxplot_all_good}, a result consistent with Figure \ref{fig:ellipse_all_good}, where all mean coefficient estimates fall within the ellipse defined in Theorem \ref{prop:cond_est_coef}.

As shown in Figure \ref{fig:boxplot_all_good}, the optimal coefficients $(\alpha_o, \beta_o)$ achieve greater variance reduction than alternative choices, with an additional 7.4\% of RVR. This highlights the advantage of jointly optimizing both coefficients for the ratio, rather than optimizing them separately without accounting for the ratio as a whole. 

For this scenario, the RAE is $2.5\%$, confirming that the $\mathcal{O}_{\infty}(n^{-3/2})$ variance approximation remains accurate even at this sample size $n=100$.

\subsection{Variance increase with non-optimal coefficients} \label{sec:worst_gordon}
The presented results are the worst-case scenario for the coefficients of \cite{gordon_efficient_1982}, happening with the covariance matrix $$ \Sigma = 
\begin{pmatrix}
1 & & & \\
0.99 & 1 & & \\
0.99 & 0.99 & 1 & \\
0.01 & 0.01 & 0.01 & 1
\end{pmatrix}.
$$

\begin{figure}[!h]
    \centering
    \begin{subfigure}[t]{0.49\linewidth}
        \centering
        \includegraphics[width=\linewidth]{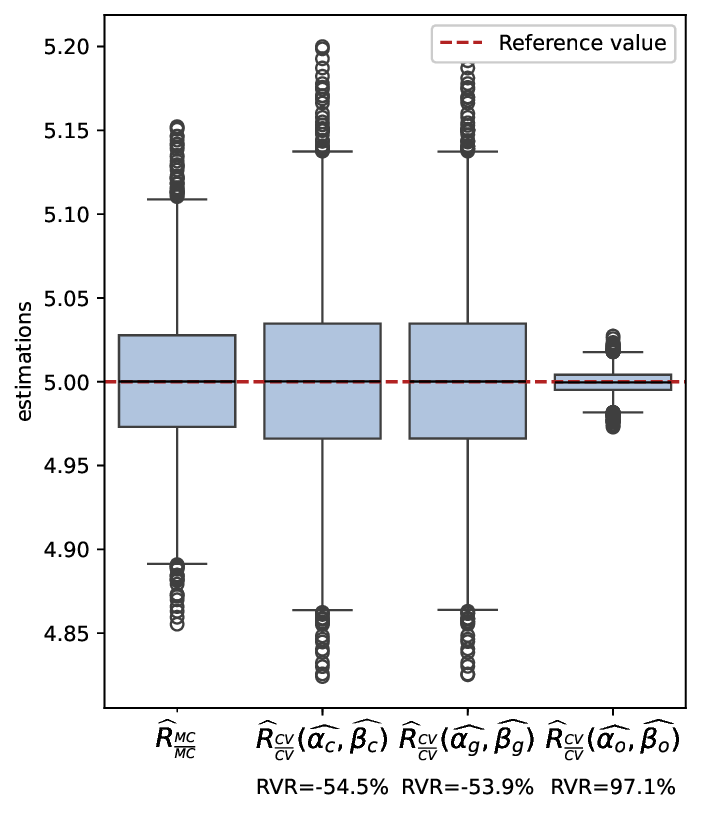}
        \subcaption{Boxplots of ratio estimations.}
        \label{fig:boxplot_negative_rvr}
    \end{subfigure}
    \hfill
    \begin{subfigure}[t]{0.49\linewidth}
        \centering
        \includegraphics[width=\linewidth]{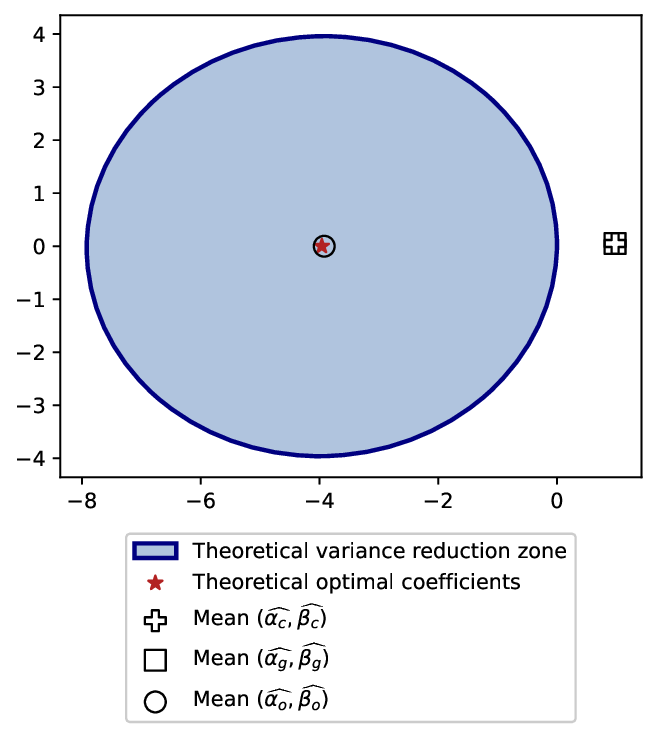}
        \subcaption{Stability ellipse and mean coefficient estimates (Theorem \ref{prop:cond_est_coef}).}
        \label{fig:ellipse_worst_gordon}
    \end{subfigure}
    \caption{\textbf{Worst-case scenario for \cite{gordon_efficient_1982} coefficients.} Only the optimal coefficients achieve variance reduction.}
\end{figure}

While the variance reduction achieved with the optimal coefficients may not always exceed that of the other coefficients, their key advantage lies in their theoretical guarantee of variance reduction. The classical coefficients and the ones from \cite{gordon_efficient_1982} can lead to a variance increase, as illustrated in Figure \ref{fig:boxplot_negative_rvr}. This result is consistent with Figure \ref{fig:ellipse_worst_gordon}, which shows that only the mean optimal coefficient estimates fall within the ellipse defined in Theorem \ref{prop:cond_est_coef}. This result highlights the importance of jointly optimizing the coefficients to minimize the variance of the whole estimator.

The $\mathcal{O}_{\infty}(n^{-3/2})$ approximation of the variance is accurate in this setting, with a measured RAE of $2.4\%$.

\subsection{Best variance reduction with the coefficients of \cite{gordon_efficient_1982}}
The presented results are the best-case scenario for the coefficients of \cite{gordon_efficient_1982}, with the following covariance matrix $$ \Sigma = 
\begin{pmatrix}
1 & & & \\
-0.71 & 1 & & \\
-0.99 & 0.67 & 1 & \\
-0.98 & 0.61 & 0.99 & 1
\end{pmatrix}.
$$

\begin{figure}[!h]
    \centering
    \begin{subfigure}[t]{0.49\linewidth}
        \centering
        \includegraphics[width=\linewidth]{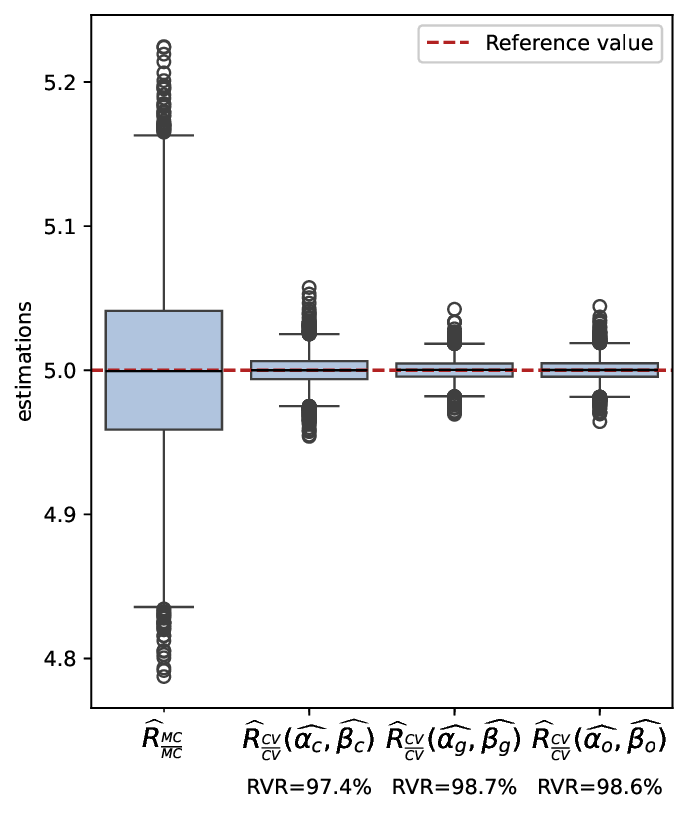}
        \subcaption{Boxplots of ratio estimations.}
        \label{fig:boxplot_gordon_best}
    \end{subfigure}
    \hfill
    \begin{subfigure}[t]{0.49\linewidth}
        \centering
        \includegraphics[width=\linewidth]{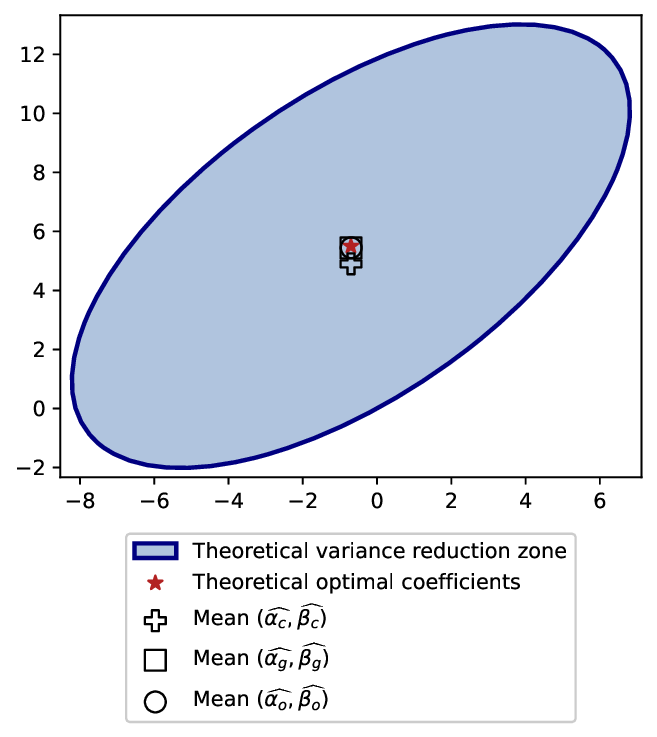}
        \subcaption{Stability ellipse and mean coefficient estimates (Theorem \ref{prop:cond_est_coef}).}
        \label{fig:ellipse_best_gordon}
    \end{subfigure}
    \caption{\textbf{Best-case scenario for \cite{gordon_efficient_1982} coefficients.} All coefficients achieve very similar variance reduction performances.}
\end{figure}

The optimal coefficients introduced in this work come with a theoretical guarantee of variance reduction. However, their complex expressions can make them harder to estimate when sample size is small (here $n=100$), which can slightly degrade empirical performance and allow the simpler \cite{gordon_efficient_1982} estimated coefficients to occasionally get closer to the true optimal coefficients.\\ Yet, even with the covariance structure optimized to favor \cite{gordon_efficient_1982} coefficients, they only outperform the optimal coefficients by a negligible 0.1\% in RVR, as illustrated in Figure \ref{fig:boxplot_gordon_best}. When \cite{gordon_efficient_1982} coefficients do successfully reduce variance, their performance is practically identical to the optimal ones, which remain robust despite estimation constraints.

Figure \ref{fig:ellipse_best_gordon} illustrates that the mean estimated coefficients are within the variance reduction ellipse. The classical coefficients remain the most peripheral, which aligns with their comparatively lower variance reduction, specifically a 1\% lower RVR.

The measured RAE is of $2.5\%$ which suggests the approximation remains accurate.

\section{Applications on real data} \label{sec:appli}
\subsection{Estimation of a proportion in a multi-fidelity aircraft design use case}
The ACV/ACV estimator is applied to estimate the strut mass fraction, defined as the ratio between the strut mass and the total mass. The strut is an external support between the fuselage and the wing that helps carry aerodynamic loads. The strut mass fraction is a quantity of interest to be optimized as it ensures the added support does not negate the aerodynamic and structural advantages of a strut-braced wing. The studied dataset contains 1,252 samples from an aircraft design study, based on \cite{meheut_strut_braced_2024}, which describes the multidisciplinary and multi-fidelity design process. The dataset combines outputs from high fidelity (HF) \cite{lannoo_mass_2025} and low fidelity (LF) \cite{carrier_multidisciplinary_2022} models to efficiently explore the design space. By applying the ACV/ACV method, the strut mass fraction can be estimated more precisely with fewer HF model evaluations, which take around 2 minutes each, leveraging the LF model, which runs in only a few seconds and serves as a control variate as in MFMC. 

Following the notation previously used in this paper, the strut mass fraction can be written as the ratio $R = \E[A]/\E[C]$ where $A$ denotes the HF strut mass and $C$ the HF total mass. The control variates $B$ and $D$ correspond to the LF strut mass and the LF total mass. Since their exact expectations are unknown, we apply the ACV method, estimating their means using the full set of $n+m=1,252$ available samples. The results are obtained by bootstrapping on 1,000 random dataset configurations of size $n=200$, and the RVR is estimated from the variance of these 1,000 estimators.

\begin{figure}[!ht]
    \centering
    \includegraphics[width=\linewidth]{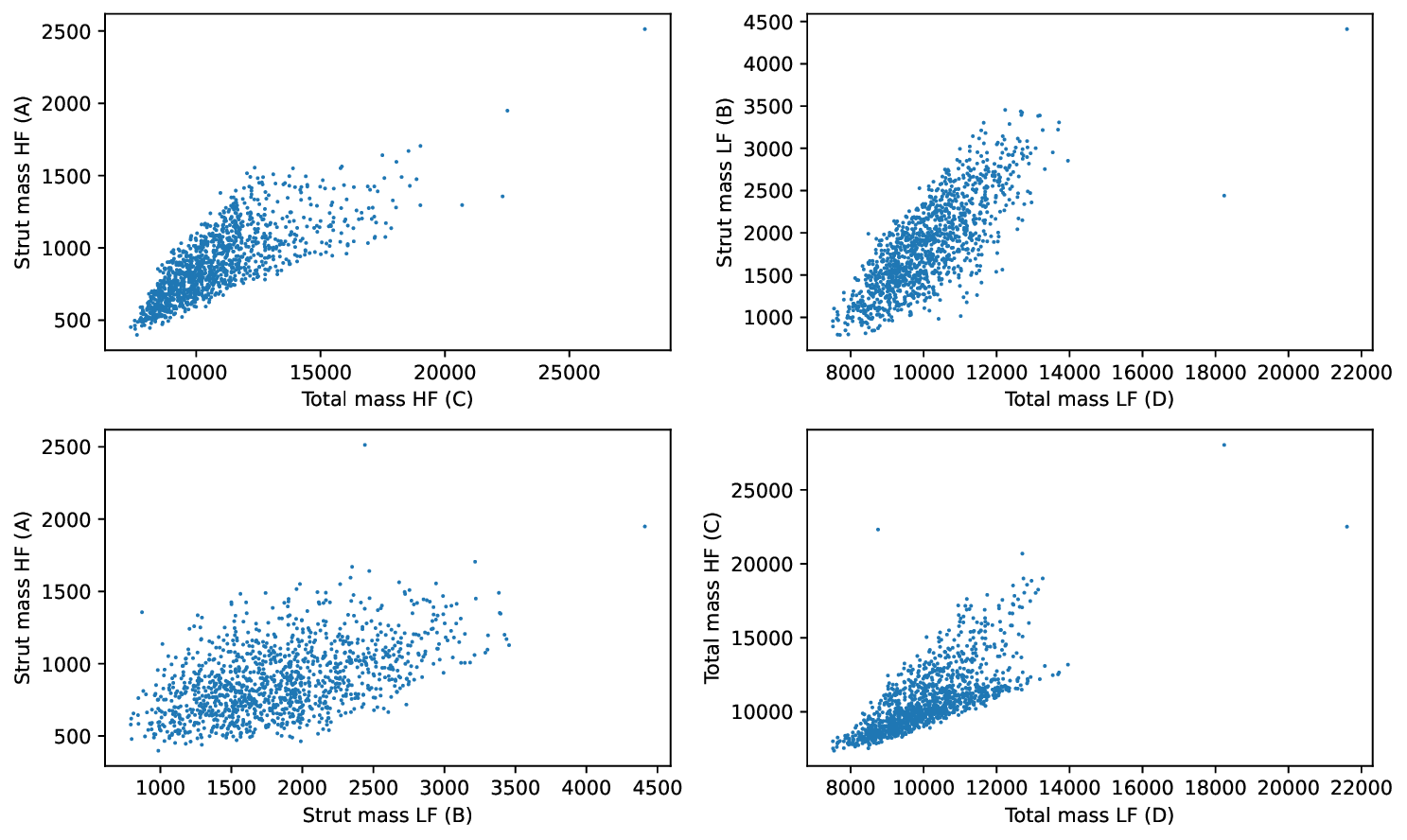}
    \vspace{-25pt}
    \caption{\textbf{Scatter plots of the 1,252 coupled aircraft design samples.} The top row illustrates the correlations between the strut mass (numerator) and total mass (denominator) within each fidelity level. The bottom row illustrates the correlation between HF variables and their respective LF control variates.}
    \label{fig:data_appli}
\end{figure}

Simulations highlighted the central role of the covariance structure between the variables and the control variates. Illustrating these dependencies, the pair plots in Figure \ref{fig:data_appli} reveal that the control variates $B$ and $D$ exhibit strong positive correlations with $A$ and $C$, suggesting a high potential for variance reduction with the CV/CV estimator.

\begin{figure}[!h]
    \centering
    \begin{subfigure}[t]{0.49\linewidth}
        \centering
        \includegraphics[width=\linewidth]{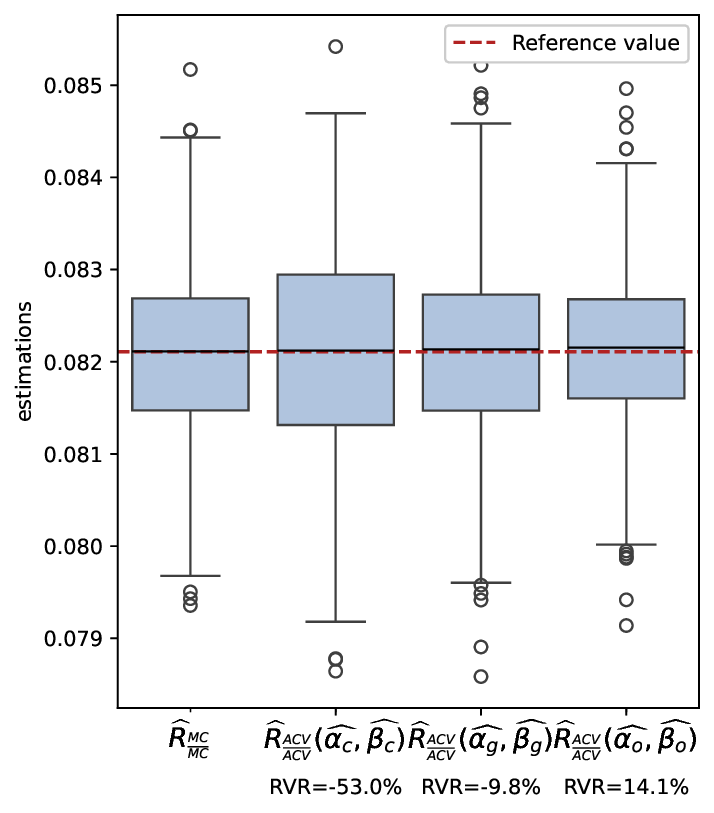}
        \subcaption{Boxplots of strut mass fraction estimations.}
        \label{fig:appli}
    \end{subfigure}
    \hfill
    \begin{subfigure}[t]{0.49\linewidth}
        \centering
        \includegraphics[width=\linewidth]{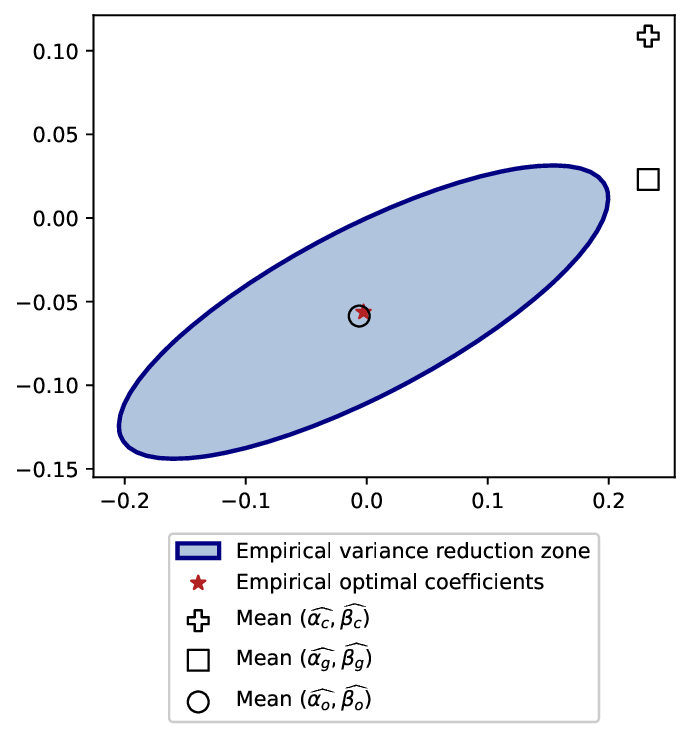}
        \subcaption{Stability ellipse and mean coefficient estimates (Theorem \ref{prop:cond_est_coef}).}
        \label{fig:ellipse_appli}
    \end{subfigure}
    \caption{\textbf{Comparative performance of the ACV/ACV coefficients for strut mass fraction estimation on the aircraft design dataset.} Only the optimal coefficients achieve variance reduction. \textit{(Reference value, variance reduction zone and optimal coefficients estimated on all 1252 available HF samples.)}}
\end{figure}

As seen in Figure \ref{fig:appli}, the variance is increased when using the non-optimal coefficients, which highlights the benefit of the proposed optimal coefficients that achieve an RVR of about 14\%. This confirms the effectiveness of the ACV/ACV estimator in improving estimation accuracy and delivering the precision of a larger HF dataset at no extra computational cost.

Figure \ref{fig:ellipse_appli} demonstrates why non-optimal coefficients lead to a variance increase, as their mean estimates fall outside the reduction boundary defined in Theorem \ref{prop:cond_est_coef}. Notably, the \cite{gordon_efficient_1982} coefficients lie closer to the ellipse, resulting in a more mitigated variance increase compared to the classical estimators.

\clearpage
\subsection{Estimation of a Conditional Value-at-Risk in a multi-fidelity electromagnetic dataset}
A CVAR \cite{hong_monte_2014}, also referred to as superquantile or expected shortfall \cite{acerbi_coherence_2002}, measures the average of values above a given quantile. While a quantile marks a threshold, the CVAR reflects the magnitude of extreme outcomes, providing a more robust assessment of risk widely applied in finance, reliability, and risk management \cite{rockafellar_optimization_2000}.

The CVAR $q^*_{\sigma}$ of $X$ of order $\sigma \in [0,1]$ is defined as \begin{align}
    q^*_{\sigma} (X) = \E[X | X \geq q_{\sigma} (X)] =\frac{\E[X \1\{X \geq q_{\sigma} (X)\}]}{\E[\1\{X \geq q_{\sigma} (X)\}]},
\end{align}
where $X \in L^2$ is an arbitrary random variable, with cumulative distribution function $F$, and $q_{\sigma} (X) = \text{inf}\{x\in \mathbb{R} | F_{X}(x) \geq \sigma\}$ is the quantile of $X$ of order $\sigma$.

The studied multi-fidelity electromagnetic dataset \cite{RSVYFT_2026} is composed of $500$ HF samples $X^{HF}$ and $4,500$ additional LF samples $X^{LF}$. Aircraft certification for electromagnetic environments requires evaluating induced constraints, in terms of currents and voltages, at embedded equipment inputs. In the case of lightning indirect effects, the waveform rise time is a critical indicator, corresponding to the time in seconds between 10\% and 90\% of the induced current waveform.

Following the notation previously used in this paper, the CVAR estimator can be written as the ratio $R = \E[A]/\E[C]$ where \begin{align*}A=X^{HF} \1\{X^{HF} \geq \widehat{q}_{\sigma} (X^{HF})\} \; \text{and} \; C=\1\{X^{HF} \geq \widehat{q}_{\sigma} (X^{HF})\}. \end{align*} Although $E[C]=\alpha$ is known, the empirical CVAR estimator is naturally expressed as a ratio of empirical means due to the estimation of the quantile. This ratio structure allows us to exploit joint fluctuations through control variates.\\
The control variates are defined analogously using the LF samples: \begin{align*}B=X^{LF} \1\{X^{LF} \geq \widehat{q}_{\sigma} (X^{LF})\} \; \text{and} \; D=\1\{X^{LF} \geq \widehat{q}_{\sigma-0.05} (X^{LF})\}. \end{align*} Since their exact expectations are unknown, we apply the ACV method, estimating their means using the full set of $n+m=5,000$ available samples. The control variate $D$ is defined with a lower quantile order $\sigma-0.05$ instead of mirroring the quantile order $\sigma$ on the target variable $C$, to avoid linear correlation between the control variates $B$ and $D$.

The results are obtained by bootstrapping on 1,000 random dataset configurations of size $n=50$, and the RVR is estimated from the variance of these 1,000 estimators.

\begin{figure}[!ht]
    \centering
    \includegraphics[height=0.2\textheight]{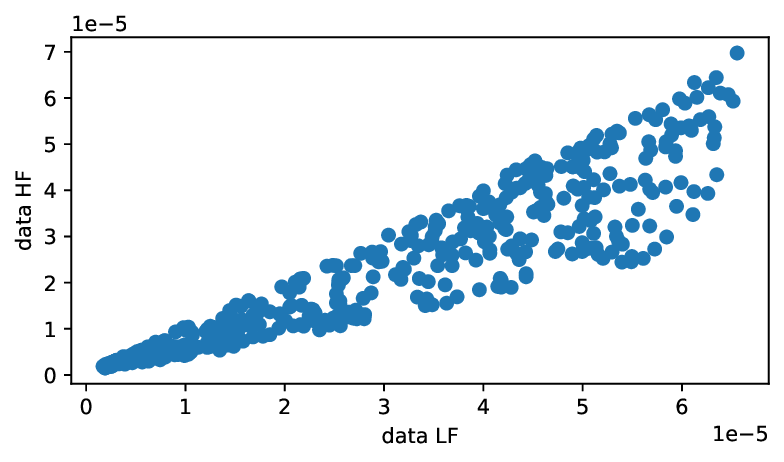}
    \vspace{-10pt}
    \caption{\textbf{Scatter plot of the 500 coupled electromagnetic samples.}}
    \label{fig:data_elec}
\end{figure}

\begin{figure}[!h]
\vspace{-40pt}
    \centering
    \begin{subfigure}[t]{0.49\linewidth}
        \centering
        \includegraphics[width=\linewidth]{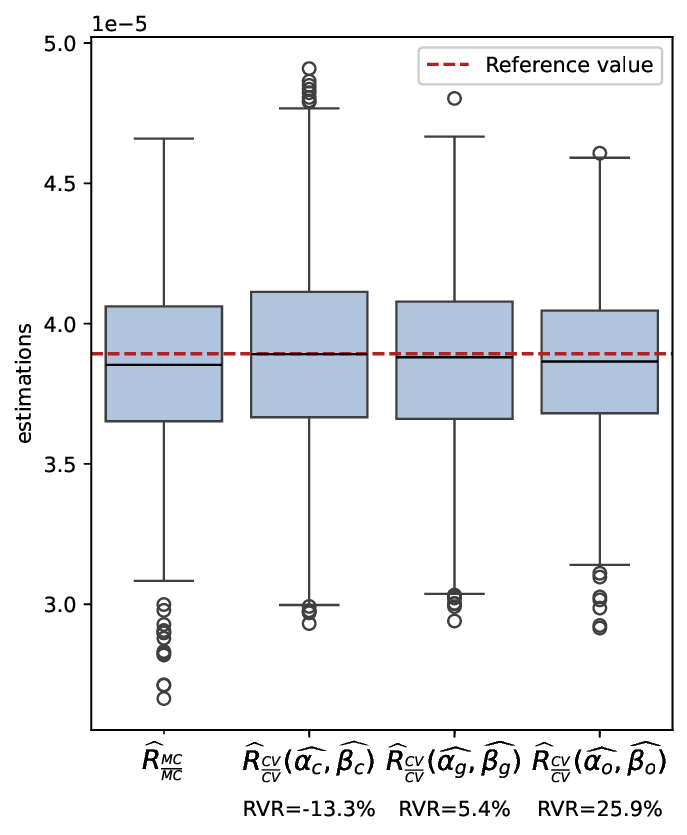}
        \subcaption{Boxplots of CVAR estimations with $\sigma=0.6$.}
        \label{fig:elec_0.6}
    \end{subfigure}
    \hfill
    \begin{subfigure}[t]{0.49\linewidth}
        \centering
        \includegraphics[width=\linewidth]{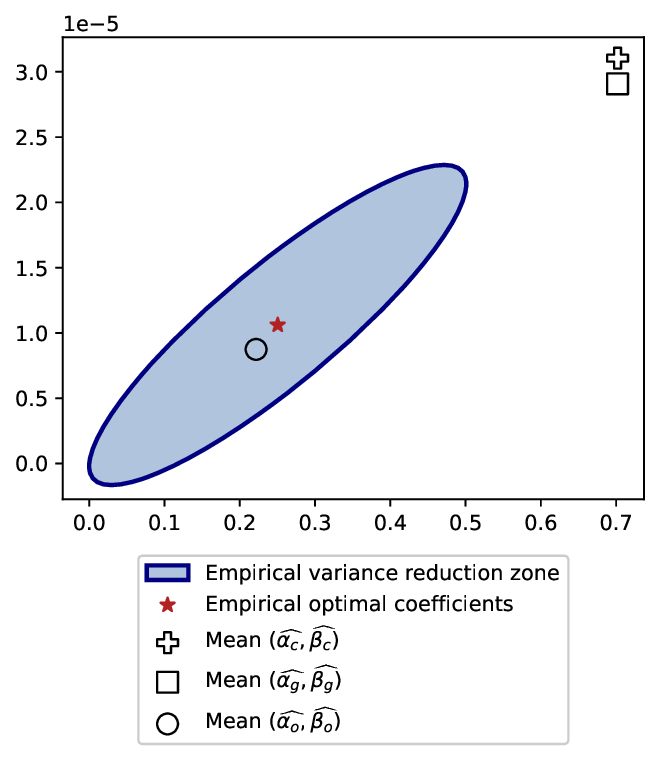}
        \subcaption{Stability ellipse and mean coefficient estimates (Theorem \ref{prop:cond_est_coef}) with $\sigma=0.6$.}
        \label{fig:ellipse_elec_0.6}
    \end{subfigure}
        \begin{subfigure}[t]{0.49\linewidth}
        \centering
        \includegraphics[width=\linewidth]{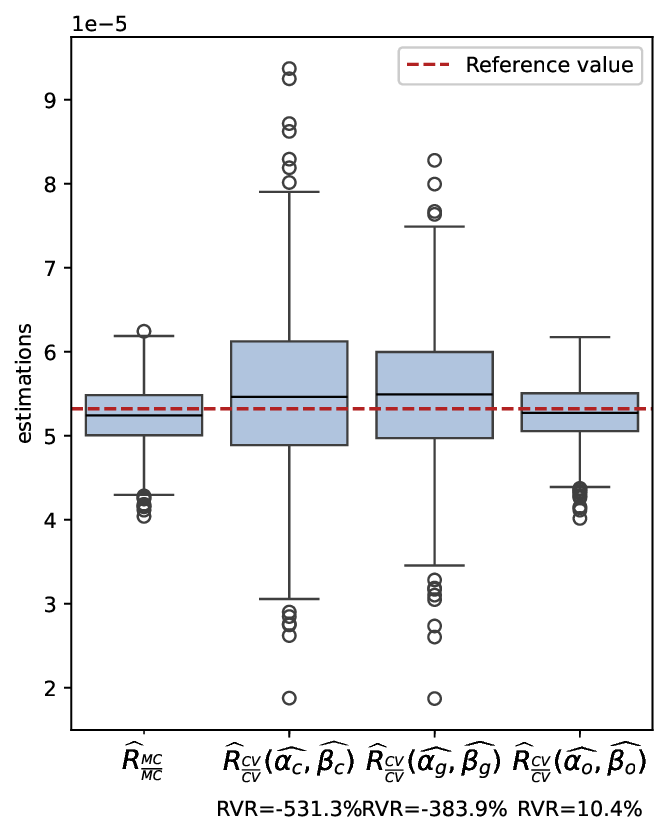}
        \subcaption{Boxplots of CVAR estimations with $\sigma=0.9$.}
        \label{fig:elec_0.9}
    \end{subfigure}
    \hfill
    \begin{subfigure}[t]{0.49\linewidth}
        \centering
        \includegraphics[width=\linewidth]{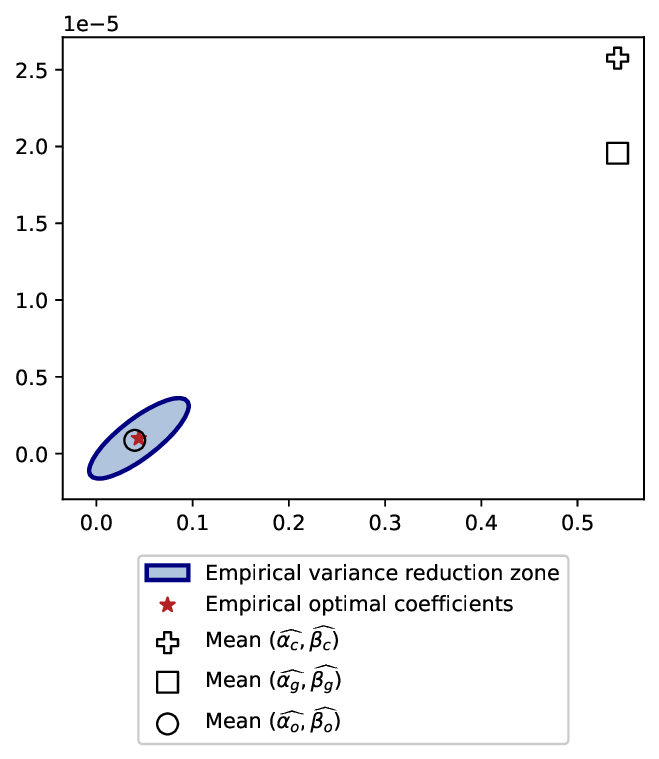}
        \subcaption{Stability ellipse and mean coefficient estimates (Theorem \ref{prop:cond_est_coef}) with $\sigma=0.9$.}
        \label{fig:ellipse_elec_0.9}
    \end{subfigure}
    \caption{\textbf{Comparative performance of the ACV/ACV coefficients for CVAR estimation on the electromagnetic dataset.} Only the optimal coefficients achieve variance reduction. \textit{(Reference value, variance reduction zone and optimal coefficients estimated on all 500 available HF samples.)}}
\end{figure}

As shown in Figure \ref{fig:data_elec}, the HF and LF samples are strongly correlated, suggesting potential variance reduction using the ACV/ACV estimator. Indeed, the estimated RVR depicted in Figures \ref{fig:elec_0.6} and \ref{fig:elec_0.9} with optimal coefficients ranges from 25\% when $\sigma=0.6$ to 10\% when $\sigma=0.9$. The variance reduction guarantee of the optimal coefficients is illustrated here, as the other sets of coefficients increase the variance.

The mean estimated non-optimal coefficients lie outside the reduction ellipse in both Figure \ref{fig:ellipse_elec_0.6} ($\sigma=0.6$) and Figure \ref{fig:ellipse_elec_0.9} ($\sigma=0.9$), which explains why they lead to a variance increase. The larger distance from the ellipse at $\sigma=0.9$ corresponds to the sharper increase in variance.

\clearpage
\section{Conclusion} 
In this paper, we introduced new optimized control variates coefficients for a variance-reduced ratio of means estimator. The resulting CV/CV estimator guarantees variance reduction when the coefficients $\alpha$ and $\beta$ are jointly optimized to minimize the overall approximated variance. Significant variance reduction can be achieved, depending on the correlation structure between the variables of interest and the control variates. In particular, we have shown that previously proposed coefficients could lead to a variance increase. To illustrate the approach, we applied the estimators with approximate control variates to estimate the strut mass fraction on a multi-fidelity aircraft design dataset and conditional values-at-risk on a multi-fidelity electromagnetic dataset.

Future work could aim to apply this method to reduce the variance of other conditional mean estimators, such as the tail dependence \cite{malevergne_investigating_2002} or the self-normalized importance sampling estimator \cite{bugallo_adaptive_2017}.\\
The method could also be adapted to estimators defined as functions of ratios. However, the proposed coefficients would not guarantee variance reduction; new coefficients would need to be derived for the full estimator. As we have shown, independently minimizing the variance of each moment may inflate the variance of a ratio; this logic extends to broader classes of estimators, such as functions of ratios or more general functions of moments.

\section*{Acknowledgments}
We gratefully acknowledge Eric Nguyen Van (ONERA/DTIS) for providing the aircraft design dataset and sharing his insights into the associated design considerations, and Isabelle Junqua (ONERA/DEMR) for providing the electromagnetic dataset along with her guidance on the relevant physical considerations.

\section*{Supplementary Materials}
All simulation codes used to produce the figures are available on GitHub at \href{https://github.com/LouB-N/Control-variates-for-variance-reduced-ratio-of-means-estimators}{https://github.com/LouB-N/}\\ \href{https://github.com/LouB-N/Control-variates-for-variance-reduced-ratio-of-means-estimators}{Control-variates-for-variance-reduced-ratio-of-means-estimators}.

\section*{Data Availability Statement}\label{data-availability-statement}
The aircraft design and electromagnetic datasets, along with the corresponding code to reproduce all figures are accessible on GitHub at \href{https://github.com/LouB-N/Control-variates-for-variance-reduced-ratio-of-means-estimators}{https://github.com/LouB-N/Control-variates-for-variance-reduced-ratio-of-} \\ \href{https://github.com/LouB-N/Control-variates-for-variance-reduced-ratio-of-means-estimators}{means-estimators}.


\section*{Declaration of generative AI use in the manuscript preparation process}
During the preparation of this work, the authors used ChatGPT 5 in order to improve language and readability. After using this tool, the authors reviewed and edited the content as needed and take full responsibility for the content of the published article.

\vspace{10pt}
\bibliographystyle{apalike}
\bibliography{references}

\newpage
\appendix
\section{Proofs} \label{sec:proofs}
\subsection{Approximated variance of the MC/MC estimator} \label{sec:proof_var_mc_mc}
An expression of the variance of the MC/MC estimator of a ratio is needed to compare it to the variance of the new variance-reduced estimators. Furthermore, the Delta Method linearization enables the derivation of closed-form expressions for the optimal coefficients that minimize the approximated variance of the CV/MC and CV/CV estimators.

Recall that the ratio to estimate is \begin{align}R = \frac{\E[A]}{\E[C]},\end{align} and its MC/MC estimator is defined as
\begin{align}\widehat{R}_{\frac{MC}{MC}} = \frac{\overline{A_n}}{\overline{C_n}},\end{align} where $\overline{X_n}=\frac{1}{n}\sum_{i=1}^{n}X_i$  for $X \in \{A,C\}$.

To derive the MC/MC estimator's variance expression, a few preliminary results are needed: 
\begin{lemma} \label{lem:mc_mc}
    \begin{enumerate}
        \item $\E[\overline{A_n}] = \E[A]$;
        \item $\E[\overline{C_n}] = \E[C]$;
        \item $\var(\overline{A_n}) = \frac{1}{n} \var(A)$;
        \item $\var(\overline{C_n}) = \frac{1}{n} \var(C)$;
        \item $\cov(\overline{A_n},\overline{C_n}) = \frac{1}{n} \cov(A,C)$.
    \end{enumerate}
\end{lemma}
\begin{proof}
    Points 1 to 4 are standard properties of expectation and variance of averages of i.i.d. random variables. \\  
    Point 5 can be proved as \begin{align*}\cov(\overline{A_n},\overline{C_n}) =& \frac{1}{n^2} \sum_{i=1}^n \sum_{j=1}^n \cov(A_i,C_j) \\=& \frac{1}{n^2} \sum_{i=1}^n \left( \cov(A_i,C_i) + \sum_{j\neq i} \cov(A_i,C_j) \right)\\=& \frac{1}{n} \cov(A,C).\end{align*}
\end{proof}

\begin{proposition}[Variance of the MC/MC estimator] \label{prop:var_mc_mc} Let $(A_i, C_i)_{i=1 \dots n}$ be i.i.d. samples from the joint distribution of the random variables $A, C \in L^2$. Assume $\E[C]\neq0$. The variance of the MC/MC estimator can be approximated as
\begin{align}
\boxed{
\vard\left(\widehat{R}_{\frac{MC}{MC}}\right) =  \frac{1}{n\E[C]^2} \left( \var(A) + R^2 \var(C) -2 RCov(A, C) \right) + \mathcal{O}_{\infty}(n^{-3/2})
}. \label{eq:var_mc_mc}
\end{align}
\end{proposition}

\begin{proof} 
Using the expressions in Lemma \ref{lem:mc_mc} and the Delta method \cite[example 5.5.27]{casella_properties_2024} \cite[section 4a.7]{sukhatme1970sampling} \cite[section 6.9]{cochran_sampling_1977}, the variance of the MC/MC estimator $\frac{\overline{A_n}}{\overline{C_n}}$ can be approximated as follows. \\
By the Central Limit Theorem, \begin{align*}
    \sqrt{n}\big((\overline{A_n}, \overline{C_n})^T - (\E[A],\E[C])^T\big) \overunderset{d}{n \to \infty}{\longrightarrow} \mathcal{N}\big(0, \Sigma \big),
\end{align*}
where $\Sigma$ is the covariance matrix of $(A,C)^T$. We apply the Delta method to find an approximation of $\var\left(\frac{\overline{A_n}}{\overline{C_n}}\right)$ as $n$ goes to infinity. \\
We define $$g(x, y) = \frac{x}{y}$$ that is twice derivable on $\mathbb{R} \times (\mathbb{R} \setminus \{0\})$. Its derivatives are: \begin{align*}
    &\frac{\partial g}{\partial x} (x,y) = \frac{1}{y}\; ; \qquad \frac{\partial g}{\partial y} (x,y) = -\frac{x}{y^2}\; ;\\
    &\frac{\partial^2 g}{\partial x^2} (x,y) = 0\; ; \qquad \frac{\partial^2 g}{\partial x \partial y} (x,y) = -\frac{1}{y^2}\; ; \qquad \frac{\partial^2 g}{\partial y^2} (x,y) = \frac{2x}{y^3} \; .
\end{align*}
Let $\Delta A = (\overline{A_n} - \E[A])$ and $\Delta C = (\overline{C_n} - \E[C])$.\\
Expanding $g(\overline{A_n}, \overline{C_n})$ in a Taylor series around $(\E[A], \E[C])$, we get: \begin{align*}
    g(\overline{A_n}, \overline{C_n}) = &g(\E[A], \E[C]) + g_x(\E[A],\E[C]) \Delta A + g_y(\E[A],\E[C]) \Delta C \\
    &+ \frac{1}{2} \left( g_{xx}(\E[A],\E[C]) (\Delta A)^2 + 2g_{xy}(\E[A],\E[C]) \Delta A \Delta C + g_{yy}(\E[A],\E[C]) (\Delta C)^2 \right)\\
    &+ \mathcal{O}\left(|\Delta A|^3 + |\Delta C|^3\right)\\
    = & R + \frac{1}{\E[C]} \Delta A -\frac{R}{\E[C]} \Delta C -\frac{1}{\E[C]^2} \Delta A \Delta C + \frac{R}{\E[C]^2} (\Delta C)^2 + \mathcal{O}_p(n^{-3/2}).
\end{align*}
Since both $\Delta A$ and $\Delta C$ are $\mathcal{O}_p(n^{-1/2})$, we have that $O\left(|\Delta A|^3 + |\Delta C|^3\right) = \mathcal{O}_p(n^{-3/2})$. \\
On one hand, we have:
\begin{align*}
    \E\big[g(\overline{A_n}, \overline{C_n})\big] = & \frac{1}{\E[C]} \E[\Delta A] -\frac{R}{\E[C]} \E[\Delta C] -\frac{1}{\E[C]^2} \E[\Delta A \Delta C] + \frac{R}{\E[C]^2} \E[(\Delta C)^2] + \mathcal{O}_{\infty}(n^{-3/2})\\
    = & R - \frac{1}{n} \frac{\cov(A,C)}{\E[C]^2} + \frac{1}{n} \frac{R}{\E[C]^2} \var(C) + \mathcal{O}_{\infty}(n^{-3/2}).
\end{align*}
Consequently, we have:
\begin{align}
    \E\big[g(\overline{A_n}, \overline{C_n})\big]^2 = & R^2 - \frac{1}{n} \frac{2R \cov(A,C)}{\E[C]^2} + \frac{1}{n} \frac{2R^2}{\E[C]^2} \var(C) + \mathcal{O}_{\infty}(n^{-3/2}). \label{eq:squared_E}
\end{align}
Since $E[\Delta A] = E[\Delta C] = 0$, the first-order terms disappear. \\
On the other hand, we have:
\begin{align*}
    g(\overline{A_n}, \overline{C_n})^2 = & R^2 + \frac{(\Delta A)^2}{\E[C]^2} + \frac{2R \Delta A}{\E[C]} - \frac{2R^2 \Delta C}{\E[C]} - \frac{4R \Delta A \Delta C}{\E[C]^2}  + \frac{3R^2 (\Delta C)^2}{\E[C]^2}  + \mathcal{O}_p(n^{-3/2}).
\end{align*}
Consequently, we have:
\begin{align}
    \E\big[g(\overline{A_n}, \overline{C_n})^2\big] = & R^2 + \frac{1}{n}\frac{\var(A)}{\E[C]^2} - \frac{1}{n}\frac{4R\cov(A,C)}{\E[C]^2} + \frac{1}{n}\frac{3R^2 \var(C)}{\E[C]^2} + \mathcal{O}_{\infty}(n^{-3/2}). \label{eq:E_squared}
\end{align}
Then gathering \eqref{eq:squared_E} and \eqref{eq:E_squared}, one can deduce that:
\begin{align}
    \vard\big(g(\overline{A_n}, \overline{C_n})^2\big) = & \E\big[g(\overline{A_n}, \overline{C_n})^2\big] - \E\big[g(\overline{A_n}, \overline{C_n})\big]^2 \nonumber\\ 
    = & \frac{1}{n \E[C]^2} \left( \var(A) - 2R\cov(A,C) + R^2 \var(C)\right) + \mathcal{O}_{\infty}(n^{-3/2}). \label{eq:var_approx}
\end{align}
The error term in \eqref{eq:var_approx} can be characterized at the next order $\mathcal{O}_{\infty}(n^{-5/2})$ as: 
\begin{align*}
    \vard\big(g(\overline{A_n}, \overline{C_n})^2\big) = & \frac{1}{n \E[C]^2} \bigg(\var(A) - 2R\cov(A,C) + R^2 \var(C)\bigg) \\
    & + \frac{1}{n^2 \E[C]^3} \begin{aligned}[t] 
        \bigg( &3R\E\big[(A-\E[A])(C-\E[C])^2\big]  -2R^2\E\big[(C-\E[C])^3\big] \\
       & -2R\E\big[(A-\E[A])^2(C-\E[C])\big]  \bigg)
    \end{aligned}[t] \\
    & - \frac{1}{n^2 \E[C]^4} \bigg( \cov(A,C) -R\var(C)\bigg)^2 + \mathcal{O}_{\infty}(n^{-5/2}). 
\end{align*}
This second-order expansion precisely quantifies the $\mathcal{O}_{\infty}(n^{-3/2})$ error of the first-order approximation, exposing two main flaws: sensitivity to the skewness of $A$ and $C$ (via third-order central moments) and a systematic overestimation bias from the variance-covariance mismatch $\big(\cov(A,C) - R\var(C)\big)^2$. In small-sample or highly asymmetric regimes, ignoring these $\mathcal{O}(n^{-2})$ dynamics makes the first-order approximation \eqref{eq:var_mc_mc} unreliable.
\end{proof}

\subsection{Approximated variance of the CV/CV estimator} \label{sec:var_cv_cv}
Recall that the CV/CV estimator is defined as \begin{align}
    \widehat{R}_{\frac{CV}{CV}}\left(\alpha,\beta\right) = \frac{\widehat{A}_{CV}(\alpha)}{\widehat{C}_{CV}(\beta)} = \frac{\overline{A_n} + \alpha (\E[B] - \overline{B_n})}{\overline{C_n} + \beta (\E[D] - \overline{D_n})}.
\end{align}

An expression of the CV/CV estimator's variance is derived, so that the optimal coefficients minimizing it can be found. To derive the CV/CV estimator's variance expression, a few preliminary results are needed:
\begin{lemma} \label{lem:cv_cv}
    \begin{enumerate}
        \item $\var(\widehat{A}_{CV}(\alpha)) = \var(\overline{A_n}) + \alpha^2 \var(\overline{B_n}) -2\alpha \cov(\overline{A_n}, \overline{B_n})$;
        \item $\var(\widehat{C}_{CV}(\beta)) = \var(\overline{C_n}) + \beta^2 \var(\overline{D_n}) -2\beta \cov(\overline{C_n}, \overline{D_n})$;
        \item $ \begin{aligned}[t]\cov(\widehat{A}_{CV}(\alpha), \widehat{C}_{CV}(\beta)) =& \cov(\overline{A_n}, \overline{C_n}) - \alpha \cov(\overline{B_n},\overline{C_n}) - \beta \cov(\overline{A_n},\overline{D_n}) + \alpha \beta \cov(\overline{B_n},\overline{D_n}). \end{aligned}$
    \end{enumerate}
\end{lemma}
\begin{proof}
    Points 1 to 3 are standard properties of variance and covariance.
\end{proof}

\begin{proposition}[Variance of the CV/CV estimator]
Let $(A_i, B_i, C_i, D_i)_{i=1 \dots n}$ be i.i.d. samples from the joint distribution of the random variables $A, B, C, D \in L^2$, where $B$ and $D$ are control variates. The variance of the CV/CV estimator can be approximated as
\begin{equation}
    \boxed{
    \begin{aligned}
    \vard\left(\widehat{R}_{\frac{CV}{CV}}\left(\alpha,\beta\right)\right) = &\frac{1}{n\E[C]^2}
    \left.
    \begin{aligned}[t]
    \Big( &\var(A) + \alpha^2\var(B) - 2\alpha\cov(A,B) + R^2\var(C)\\
    &+ \beta^2R^2\var(D) - 2\beta R^2\cov(C,D)  - 2R\cov(A,C) \\
    &+ 2\alpha R \cov(B,C)+ 2\beta R \cov(A,D) - 2\alpha\beta R\cov(B,D) \Big)
    \end{aligned}
    \right. \\
    &+ \mathcal{O}_{\infty}(n^{-3/2}).
    \end{aligned}
    } \label{eq:var_cv_cv}
\end{equation}
\end{proposition}

\begin{proof}
    Using the results from Lemma \ref{lem:cv_cv} and Proposition \ref{prop:var_mc_mc}, the CV/CV estimator's variance is derived as:
    \begin{align*}
        \vard\left(\widehat{R}_{\frac{CV}{CV}}\left(\alpha,\beta\right)\right)
        =&  \frac{1}{\E[\widehat{C}_{CV}(\beta)]^2}\var(\widehat{A}_{CV}(\alpha)) + \frac{\E[\widehat{A}_{CV}(\alpha)]^2}{\E[\widehat{C}_{CV}(\beta)]^4} \var(\widehat{C}_{CV}(\beta)) \\
        & -2 \frac{\E[\widehat{A}_{CV}(\alpha)]}{\E[\widehat{C}_{CV}(\beta)]^3} \cov\left(\widehat{A}_{CV}(\alpha),\widehat{C}_{CV}(\beta)\right) + \mathcal{O}_{\infty}(n^{-3/2})\\
        =& \frac{1}{\E[C]^2}\var(\widehat{A}_{CV}(\alpha)) + \frac{\E[A]^2}{\E[C]^4} \var(\widehat{C}_{CV}(\beta)) \\
        &-2 \frac{\E[A]}{\E[C]^3} \cov\left(\widehat{A}_{CV}(\alpha),\widehat{C}_{CV}(\beta)\right) + \mathcal{O}_{\infty}(n^{-3/2})\\
        =& \frac{1}{\E[C]^2} \var(\overline{A_n}) + \frac{\alpha^2}{\E[C]^2} \var(\overline{B_n}) - \frac{2\alpha}{\E[C]^2} \cov(\overline{A_n}, \overline{B_n}) + \frac{\E[A]^2}{\E[C]^4} \var(\overline{C_n}) \\
        & + \beta^2 \frac{\E[A]^2}{\E[C]^4} \var(\overline{D_n}) -2\beta \frac{\E[A]^2}{\E[C]^4} \cov(\overline{C_n}, \overline{D_n}) -2 \frac{\E[A]}{\E[C]^3} \cov(\overline{A_n}, \overline{C_n}) \\
        & +2 \alpha \frac{\E[A]}{\E[C]^3}\cov(\overline{B_n},\overline{C_n}) +2 \beta \frac{\E[A]}{\E[C]^3}\cov(\overline{A_n},\overline{D_n})\\
        &-2\alpha \beta \frac{\E[A]}{\E[C]^3}\cov(\overline{B_n},\overline{D_n}) + \mathcal{O}_{\infty}(n^{-3/2})
    \end{align*}
    \begin{align*}
        \vard\left(\widehat{R}_{\frac{CV}{CV}}\left(\alpha,\beta\right)\right)
        =& \frac{1}{\E[C]^2} \frac{1}{n}\var(A) + \frac{\alpha^2}{\E[C]^2} \frac{1}{n}\var(B) - \frac{2\alpha}{\E[C]^2} \frac{1}{n}\cov(A, B) + \frac{R^2}{\E[C]^2} \frac{1}{n}\var(C) \\
        & + \beta^2 \frac{R^2}{\E[C]^2} \frac{1}{n}\var(D) -2\beta \frac{R^2}{\E[C]^2} \frac{1}{n}\cov(C, D) -2 \frac{R}{\E[C]^2} \frac{1}{n}\cov(A, C) \\
        & +2 \alpha \frac{R}{\E[C]^2} \frac{1}{n}\cov(B,C) +2 \beta \frac{R}{\E[C]^2} \frac{1}{n}\cov(A,D)\\
        &-2\alpha \beta \frac{R}{\E[C]^2} \frac{1}{n}\cov(B,D) + \mathcal{O}_{\infty}(n^{-3/2}).
    \end{align*}
\end{proof}

\vspace{10pt}
\subsection{Optimal coefficients of the CV/CV estimator} \label{sec:cv_cv_opt_coef}
\begin{proof}[Proof of Proposition \ref{prop_cv_cv}]
We recall the definition of the optimal coefficients:
\begin{align*}
    (\alpha_o, \beta_o) := \underset{(\alpha,\beta) \in \mathbb{R}^2}{\text{argmin}} \; \vard\left(\frac{\overline{A_n} + \alpha (\E[B] - \overline{B_n})}{\overline{C_n} + \beta (\E[D] - \overline{D_n})}\right).
    \end{align*}
The approximated variance in \eqref{eq:var_cv_cv} is minimized by differentiating its expression with respect to $\alpha$ and $\beta$, and solving for the optimal value:
\begin{align*}
    &\nabla \vard\left(\widehat{R}_{\frac{CV}{CV}}\left(\alpha,\beta\right)\right) \\
    =& \begin{bmatrix}
    \frac{\partial }{\partial \alpha} \vard\left(\widehat{R}_{\frac{CV}{CV}}\left(\alpha,\beta\right)\right)\\
    \frac{\partial }{\partial \beta} \vard\left(\widehat{R}_{\frac{CV}{CV}}\left(\alpha,\beta\right)\right)
    \end{bmatrix} \\
    =& \begin{bmatrix}
    \frac{2 \alpha}{\E[C]^2} \var(\overline{B_n}) - \frac{2}{\E[C]^2} \cov(\overline{A_n}, \overline{B_n}) +2 \frac{\E[A]}{\E[C]^3}\cov(\overline{B_n},\overline{C_n}) -2\beta \frac{\E[A]}{\E[C]^3}\cov(\overline{B_n},\overline{D_n})\\
    2 \beta \frac{\E[A]^2}{\E[C]^4} \var(\overline{D_n}) -2 \frac{\E[A]^2}{\E[C]^4} \cov(\overline{C_n}, \overline{D_n}) +2 \frac{\E[A]}{\E[C]^3}\cov(\overline{A_n},\overline{D_n}) - 2\alpha \frac{\E[A]}{\E[C]^3}\cov(\overline{B_n},\overline{D_n})
    \end{bmatrix}.
\end{align*}

The Jacobian of $\vard\left(\widehat{R}_{\frac{CV}{CV}}\left(\alpha,\beta\right)\right)$ with respect to $(\alpha,\beta)$ is given by:
\begin{align*}
    \nabla^2 \vard\left(\widehat{R}_{\frac{CV}{CV}}\left(\alpha,\beta\right)\right) &= 
    \begin{bmatrix}
        \frac{\partial^2 }{\partial \alpha^2} \vard\left(\widehat{R}_{\frac{CV}{CV}}\left(\alpha,\beta\right)\right) & \frac{\partial^2 }{\partial \alpha \partial \beta} \vard\left(\widehat{R}_{\frac{CV}{CV}}\left(\alpha,\beta\right)\right) \\
        \frac{\partial^2 }{\partial \beta \partial \alpha} \vard\left(\widehat{R}_{\frac{CV}{CV}}\left(\alpha,\beta\right)\right) & \frac{\partial^2 }{\partial \beta^2} \vard\left(\widehat{R}_{\frac{CV}{CV}}\left(\alpha,\beta\right)\right)
    \end{bmatrix} \\
    &= \begin{bmatrix}
        \frac{2}{\E[C]^2} \var(\overline{B_n}) & -2 \frac{\E[A]}{\E[C]^3}\cov(\overline{B_n},\overline{D_n}) \\
        -2 \frac{\E[A]}{\E[C]^3}\cov(\overline{B_n},\overline{D_n}) & 2 \frac{\E[A]^2}{\E[C]^4} \var(\overline{D_n})
    \end{bmatrix}\\
    &:= M.
\end{align*}
The matrix M is positive definite as
\begin{align*}
    det (M) >0
    \Leftrightarrow& \frac{2}{\E[C]^2} \var(\overline{B_n}) \times 2 \frac{\E[A]^2}{\E[C]^4} \var(\overline{D_n}) - \left(-2 \frac{\E[A]}{\E[C]^3}\cov(\overline{B_n},\overline{D_n})\right)^2 >0\\
    \Leftrightarrow& \frac{4\E[A]^2}{\E[C]^6} \var(\overline{B_n})\var(\overline{D_n}) - \frac{4\E[A]^2}{\E[C]^6} \cov(\overline{B_n},\overline{D_n})^2 >0 \\
    \Leftrightarrow& \var(\overline{B_n})\var(\overline{D_n}) > \cov(\overline{B_n},\overline{D_n})^2 \\
    \Leftrightarrow& |\corr(B,D)|<1.
\end{align*}
Therefore, the function $(\alpha,\beta) \mapsto \vard\left(\widehat{R}_{\frac{CV}{CV}}\left(\alpha,\beta\right)\right)$ is strictly convex. So the values of $\alpha_o$ and $\beta_o$ for which $\nabla \vard\left(\widehat{R}_{\frac{CV}{CV}}\left(\alpha_o,\beta_o\right)\right) = \begin{bmatrix} 0\\0\end{bmatrix}$ are the unique minimizer of the variance, under the mild condition $|\corr(B,D)|<1$. Other settings with $|\corr(B,D)|=1$ are explored in the next section.

The optimal coefficients $\alpha_o$ and $\beta_o$ are solution of the following system:
$$
\nabla \vard\left(\widehat{R}_{\frac{CV}{CV}}\left(\alpha_o,\beta_o\right)\right) = \begin{bmatrix}
    0 \\
    0
\end{bmatrix}
\Leftrightarrow M \cdot \begin{bmatrix}
    \alpha_o \\
    \beta_o
\end{bmatrix} = \begin{bmatrix}
    c_1 \\
    c_2
\end{bmatrix},$$

with $M = \begin{bmatrix}
    \frac{2}{\E[C]^2} \var(\overline{B_n}) & -2\frac{\E[A]}{\E[C]^3} \cov(\overline{B_n}, \overline{D_n}) \\
    -2\frac{\E[A]}{\E[C]^3} \cov(\overline{B_n}, \overline{D_n}) & 2\frac{\E[A]^2}{\E[C]^4} \var(\overline{D_n})
\end{bmatrix},$

and $\begin{bmatrix}
    c_1 \\
    c_2
\end{bmatrix}
= \begin{bmatrix}
    \frac{2}{\E[C]^2} \cov(\overline{A_n}, \overline{B_n}) - 2 \frac{\E[A]}{\E[C]^3} \cov(\overline{B_n}, \overline{C_n}) \\
    2 \frac{\E[A]^2}{\E[C]^4} \cov(\overline{C_n}, \overline{D_n}) - 2 \frac{\E[A]}{\E[C]^3} \cov(\overline{A_n}, \overline{D_n})
\end{bmatrix}$.\\

The solution is given by: $
\begin{bmatrix}
    \alpha_o \\
    \beta_o
\end{bmatrix}
= M^{-1} \cdot \begin{bmatrix}
    c_1 \\
    c_2
\end{bmatrix}$, where the inverse of $M$ is:\\
$M^{-1} =\frac{1}{\frac{4\E[A]^2}{\E[C]^6} \var(\overline{B_n})\var(\overline{D_n}) - \frac{4\E[A]^2}{\E[C]^6} \cov(\overline{B_n},\overline{D_n})^2} \begin{bmatrix}
    2\frac{\E[A]^2}{\E[C]^4} \var(\overline{D_n}) & 2\frac{\E[A]}{\E[C]^3} \cov(\overline{B_n}, \overline{D_n}) \\
    2\frac{\E[A]}{\E[C]^3} \cov(\overline{B_n}, \overline{D_n}) & \frac{2}{\E[C]^2} \var(\overline{B_n})
\end{bmatrix}$. \\

The coefficients can be derived as follows: \small
\begin{flushleft}
\hspace*{-90pt}
\begin{minipage}{1.4\textwidth}
\begin{align*}
    \alpha_o &= \frac{2\frac{\E[A]^2}{\E[C]^4} \var(\overline{D_n})\left( \frac{2}{\E[C]^2} \cov(\overline{A_n}, \overline{B_n}) - 2 \frac{\E[A]}{\E[C]^3} \cov(\overline{B_n}, \overline{C_n}) \right) + 2\frac{\E[A]}{\E[C]^3} \cov(\overline{B_n}, \overline{D_n})\left( 2 \frac{\E[A]^2}{\E[C]^4} \cov(\overline{C_n}, \overline{D_n}) - 2 \frac{\E[A]}{\E[C]^3} \cov(\overline{A_n}, \overline{D_n}) \right)}{\frac{4\E[A]^2}{\E[C]^6} \var(\overline{B_n})\var(\overline{D_n}) - \frac{4\E[A]^2}{\E[C]^6} \cov(\overline{B_n},\overline{D_n})^2} \\
    &= \frac{2\frac{R^2}{\E[C]^2} \var(D)\left( \frac{2}{\E[C]^2} \cov(A, B) - 2 \frac{R}{\E[C]^2} \cov(B, C) \right) + 2\frac{R}{\E[C]^2} \cov(B, D)\left( 2 \frac{R^2}{\E[C]^2} \cov(C, D) - 2 \frac{R}{\E[C]^2} \cov(A, D) \right)}{\frac{4R^2}{\E[C]^4} \var(B)\var(D) - \frac{4R^2}{\E[C]^4} \cov(B,D)^2} \\
    &= \frac{\frac{4R^2}{\E[C]^4} \var(D)\cov(A, B) - \frac{4R^3}{\E[C]^4} \var(D)\cov(B, C)  + \frac{4R^3}{\E[C]^4} \cov(B, D)\cov(C, D) - \frac{4R^2}{\E[C]^4} \cov(B, D)\cov(A, D)}{\frac{4R^2}{\E[C]^4} \var(B)\var(D) - \frac{4R^2}{\E[C]^4} \cov(B,D)^2} \\
    &= \frac{\var(D)\cov(A, B) - R \var(D)\cov(B, C)  + R \cov(B, D)\cov(C, D) -  \cov(B, D)\cov(A, D)}{\var(B)\var(D) - \cov(B,D)^2};
\end{align*}
\end{minipage}
\end{flushleft}

\begin{flushleft}
\hspace*{-90pt}
\begin{minipage}{1.4\textwidth}
\begin{align*}
    \beta_o &= \frac{2\frac{\E[A]}{\E[C]^3} \cov(\overline{B_n}, \overline{D_n})\left( \frac{2}{\E[C]^2} \cov(\overline{A_n}, \overline{B_n}) - 2 \frac{\E[A]}{\E[C]^3} \cov(\overline{B_n}, \overline{C_n}) \right) + \frac{2}{\E[C]^2} \var(\overline{B_n})\left( 2 \frac{\E[A]^2}{\E[C]^4} \cov(\overline{C_n}, \overline{D_n}) - 2 \frac{\E[A]}{\E[C]^3} \cov(\overline{A_n}, \overline{D_n}) \right)}{\frac{4\E[A]^2}{\E[C]^6} \var(\overline{B_n})\var(\overline{D_n}) - \frac{4\E[A]^2}{\E[C]^6} \cov(\overline{B_n},\overline{D_n})^2} \\
    &= \frac{2\frac{R}{\E[C]^2} \cov(B, D)\left( \frac{2}{\E[C]^2} \cov(A, B) - 2 \frac{R}{\E[C]^2} \cov(B, C) \right) + \frac{2}{\E[C]^2} \var(B)\left( 2 \frac{R^2}{\E[C]^2} \cov(C, D) - 2 \frac{R}{\E[C]^2} \cov(A, D) \right)}{\frac{4R^2}{\E[C]^4} \var(B)\var(D) - \frac{4R^2}{\E[C]^4} \cov(B,D)^2} \\
    &= \frac{\frac{4R}{\E[C]^4} \cov(B, D)\cov(A, B) - \frac{4R^2}{\E[C]^4} \cov(B, D) \cov(B, C) + \frac{4R^2}{\E[C]^4} \var(B)\cov(C, D) - \frac{4R}{\E[C]^4} \var(B)\cov(A, D)}{\frac{4R^2}{\E[C]^4} \var(B)\var(D) - \frac{4R^2}{\E[C]^4} \cov(B,D)^2} \\
    &= \frac{\frac{1}{R} \cov(B, D)\cov(A, B) -  \cov(B, D) \cov(B, C) + \var(B)\cov(C, D) - \frac{1}{R} \var(B)\cov(A, D)}{\var(B)\var(D) - \cov(B,D)^2}.
\end{align*}
\end{minipage}
\end{flushleft}
\normalsize

The optimal coefficients are:
\begin{align}
    \boxed{\alpha_o = \frac{\var(D)\cov(A, B) - R \var(D)\cov(B, C)  + R \cov(B, D)\cov(C, D) -  \cov(B, D)\cov(A, D)}{\var(B)\var(D) - \cov(B,D)^2}}; \\
    \boxed{\beta_o = \frac{\frac{1}{R} \cov(B, D)\cov(A, B) -  \cov(B, D) \cov(B, C) + \var(B)\cov(C, D) - \frac{1}{R} \var(B)\cov(A, D)}{\var(B)\var(D) - \cov(B,D)^2}}.
\end{align}
\end{proof}

\vspace{10pt}
\subsection{Comparing the approximated variances of the MC/MC and the CV/CV estimators with different coefficients} \label{sec:diff_var}
The approximated difference between variances can be expressed as:
\begin{align} \label{eq:diff_cv_cv_mv_mv}
&\vard\left(\widehat{R}_{\frac{CV}{CV}}\left(\alpha,\beta\right)\right) - \vard\left(\widehat{R}_{\frac{MC}{MC}}\right) \nonumber \\
&= \frac{1}{n\E[C]^2}
\left(
\begin{aligned}[t]
&\alpha^2 \var(B) - 2\alpha \cov(A, B) + \beta^2 R^2 \var(D) - 2\beta R^2 \cov(C, D) \\
&\quad + 2 \alpha R \cov(B, C) + 2 \beta R \cov(A, D) - 2\alpha \beta R \cov(B, D)
\end{aligned}
\right) + \mathcal{O}_{\infty}(n^{-3/2}).
\end{align}

$\bullet$ With the classical coefficients $\alpha_c = \frac{\cov(A, B)}{\var(B)}$ and $\beta_c = \frac{\cov(C, D)}{\var(D)}$, the approximated variance difference is:
\begin{align*}
&\vard\left(\widehat{R}_{\frac{CV}{CV}}\left(\alpha_c,\beta_c\right)\right) - \vard\left(\widehat{R}_{\frac{MC}{MC}}\right) \\
= &\frac{1}{n\E[C]^2 \var(B)\var(D)}
\left.
\begin{aligned}[t]
\big( &-\cov(A,B)^2 \var(D) - R^2 \cov(C,D)^2 \var(B)\\
&+ 2R \cov(A,B) \cov(B,C) \var(D) + 2R \cov(C,D) \cov(A,D) \var(B)\\
& - 2R \cov(A,B) \cov(C,D) \cov(B,D) \; \big) + \mathcal{O}_{\infty}(n^{-3/2})
\end{aligned}
\right. \\
:= &\frac{1}{n\E[C]^2} \frac{\textbf{T}}{\var(B)\var(D)} + \mathcal{O}_{\infty}(n^{-3/2}).
\end{align*}
There is a variance reduction if $\textbf{T}<0$. \\

$\bullet$ With the coefficients from \cite{gordon_efficient_1982} $\alpha_g = \frac{\cov(A, B)}{\var(B)}$ and $\beta_g = \frac{\cov(C, D) - \frac{1}{R} \cov(A,D) + \frac{\alpha}{R} \cov(B,D)}{\var(D)}$, the approximated variance difference is:
\begin{align*}
&\vard\left(\widehat{R}_{\frac{CV}{CV}}\left(\alpha_g,\beta_g\right)\right) - \vard\left(\widehat{R}_{\frac{MC}{MC}}\right) = \frac{1}{n\E[C]^2} \begin{aligned}[t]\bigg( &\frac{\textbf{T} - \cov(A,D)^2Var(B) + 2Cov(A,B)\cov(B,D)\cov(A,D)}{\var(B)\var(D)} \\&- \frac{\cov(A,B)^2Cov(B,D)^2}{\var(B)^2Var(D)} \bigg) + \mathcal{O}_{\infty}(n^{-3/2}). \end{aligned} \end{align*}
There is a variance reduction if $\var(B) \bigg(\textbf{T} - \cov(A,D)^2Var(B) + 2Cov(A,B)\cov(B,D)\cov(A,D)\bigg) - \cov(A,B)^2Cov(B,D)^2 < 0$. \\

$\bullet$ With the optimal coefficients $\alpha_o = \frac{\var(D)\cov(A, B) - R \var(D)\cov(B, C)  + R \cov(B, D)\cov(C, D) -  \cov(B, D)\cov(A, D)}{\var(B)\var(D) - \cov(B,D)^2} \text{ and }\\ \beta_o = \frac{\frac{1}{R} \cov(B, D)\cov(A, B) -  \cov(B, D) \cov(B, C) + \var(B)\cov(C, D) - \frac{1}{R} \var(B)\cov(A, D)}{\var(B)\var(D) - \cov(B,D)^2}$, the approximated variance difference is: 
\begin{align*} &\vard\left(\widehat{R}_{\frac{CV}{CV}}\left(\alpha_o,\beta_o\right)\right) - \vard\left(\widehat{R}_{\frac{MC}{MC}}\right)\\
&= - \frac{1}{n\E[C]^2} \frac{Var\bigg(\big(R\cov(B,C)-\cov(A,B)\big)D - \big(R\cov(C,D)-\cov(A,D)\big)B\bigg)}{\var(B)\var(D)-\cov(B,D)^2} + \mathcal{O}_{\infty}(n^{-3/2})\leq 0 .\end{align*}
The variance reduction is guaranteed.

\vspace{10pt}
\subsection{Condition for variance reduction with estimated control variates coefficients} 
\subsubsection{Mean estimator} \label{sec:cond_est_coef_cv}
Let $\widehat{A}_{CV}(\alpha_c) = \overline{A_n} + \alpha_c (\E[B] - \overline{B_n})$ be the CV mean estimator of $\E[A]$ with optimal coefficient $\alpha_c = \frac{\cov(A, B)}{\var(B)}$ \eqref{eq:classic_alpha}.

The variance reduction with an estimated coefficient $\widehat{\alpha}_c$ is conditioned by:
\begin{align*}
    \var\left(\widehat{A}_{CV}(\widehat{\alpha}_c)\right) < \var\left(\overline{A_n}\right)
    \Longleftrightarrow& \; \var\left(\overline{A_n}\right) - 2 \widehat{\alpha}_c \cov\left(\overline{A_n},\overline{B_n}\right) + \widehat{\alpha}_c^2 \var\left(\overline{B_n}\right) < \var\left(\overline{A_n}\right) \\
    \Longleftrightarrow& \; - 2 \widehat{\alpha}_c \cov\left(A,B\right) + \widehat{\alpha}_c^2 \var\left(B\right) < 0 \\
    \Longleftrightarrow& \; -2\widehat{\alpha}_c \alpha_c + \widehat{\alpha}_c^2 < 0\\
    \Longleftrightarrow& \; \widehat{\alpha}_c \left(\widehat{\alpha}_c -2 \alpha_c\right) < 0\\
    \Longleftrightarrow& \;  \widehat{\alpha}_c \in (\min(0, 2\alpha_c), \max(0, 2\alpha_c)).
\end{align*}

\subsubsection{Ratio of means estimator} \label{sec:cond_est_coef_cv_cv}
This section details the proof of Theorem \ref{prop:cond_est_coef}.\\
Let $\widehat{R}_{\frac{CV}{CV}}(\alpha, \beta) = \frac{\overline{A_n} + \alpha (\E[B] - \overline{B_n})}{\overline{C_n} + \beta (\E[D] - \overline{D_n})}$ be the CV/CV ratio of means estimator of $R= \frac{\E[A]}{\E[C]}$, with optimal coefficients $\alpha_o$ and $\beta_o$ as given in \eqref{eq:alpha_opt} and \eqref{eq:beta_opt}.

The approximated variance reduction with estimated coefficients $\widehat{\alpha}$ and $\widehat{\beta}$ is conditioned by:
\begin{align}
    \vard\left(\widehat{R}_{\frac{CV}{CV}}(\widehat{\alpha}, \widehat{\beta})\right) < \vard\left(\widehat{R}_{\frac{MC}{MC}}\right) 
    \Longleftrightarrow& \; \widehat{\alpha}^2 \var(B) - 2\widehat{\alpha} \cov(A, B) + \widehat{\beta}^2 R^2 \var(D) - 2\widehat{\beta} R^2 \cov(C, D) \nonumber \\
    &\quad + 2 \widehat{\alpha} R \cov(B, C) + 2 \widehat{\beta} R \cov(A, D) - 2\widehat{\alpha} \widehat{\beta} R \cov(B, D) <0, \label{eq:var_diff}
\end{align}
where the second expression is derived with the first-order approximation of the variance of the ratio \eqref{eq:var_cv_cv} as in \eqref{eq:diff_cv_cv_mv_mv}.\\
Defining the vector of estimated coefficients $\widehat{x} = \left(\widehat{\alpha}, R\widehat{\beta}\right)^T$ and the covariance matrix \\$\Lambda = \begin{pmatrix} \var(B) & -\cov(B,D) \\ -\cov(B,D) & \var(D) \end{pmatrix}$, the inequality can be written as:
\begin{align*}
    \widehat{x}^T \Lambda \widehat{x} - 2b^T \widehat{x} <0,
\end{align*}
where $b^T = \left(\cov(A,B)-R\cov(B,C), \;R\cov(C,D)-\cov(A,D)\right)$.\\
Following the notation in Appendix \ref{sec:cv_cv_opt_coef}, we identify $\Lambda = \frac{n \E[C]^2}{2} M$ and $b = \frac{n \E[C]^2}{2} (c_1, R c_2)^T$. The optimal coefficients $\alpha_o$ and $\beta_o$ are defined by the relation $M (\alpha_o, \beta_o)^T = (c_1, c_2)^T$. Consequently, if we define the optimal vector as $x_o = (\alpha_o, R\beta_o)^T$, it satisfies the system $\Lambda x_o = b$. Substituting $b^T = x_o^T \Lambda$ into the inequality \eqref{eq:var_diff} yields:
\begin{align*}
    \widehat{x}^T \Lambda \widehat{x} - 2x_o^T \Lambda \widehat{x} <0
    \Longleftrightarrow &\; \widehat{x}^T \Lambda \widehat{x} - 2x_o^T \Lambda \widehat{x} + x_o^T \Lambda x_o < x_o^T \Lambda x_o \\
    \Longleftrightarrow& \; (\widehat{x}-x_o)^T \Lambda (\widehat{x}-x_o) < x_o^T \Lambda x_o \, .
\end{align*}

\vspace{10pt}
\subsection{Approximate control variates}  \label{sec:approx}
Recall that the exact control variates estimator of $\E[A]$ is $\widehat{A}_{CV}(\alpha) = \overline{A_n} + \alpha (\E[B] - \overline{B_n})$ and the approximate estimator is $\widehat{A}_{ACV}(\alpha) = \overline{A_n} + \alpha (\overline{B_{n+m}} - \overline{B_n})$, where $B$ is the control variate.

For $X \in \{A,C\}$ and $Y \in \{B,D\}$, we have the following identities:
\begin{align*}
&\var(\overline{Y_{n+m}} - \overline{Y_n}) = \frac{m}{n+m} \var(\E[Y] - \overline{Y_n});\\
&\cov(\overline{Y_{n+m}} - \overline{Y_n}, \overline{X_n}) = \frac{m}{n+m} \cov(\E[Y] - \overline{Y_n}, \overline{X_n});\\
&\cov(\overline{B_{n+m}} - \overline{B_n}, \overline{D_{n+m}} - \overline{D_n}) = \frac{m}{n+m} \cov(\E[B] - \overline{B_n}, \E[D] - \overline{D_n}).
\end{align*}

The coefficients that minimize the variance remain the same in the approximate case as in the exact case.

It can easily be shown that the variance reduction in the approximate case is the same as in the exact case, with an additional factor of $\frac{m}{n+m}$ i.e.
\begin{align}
\vard\left(\widehat{R}_{\frac{ACV}{ACV}}\left(\alpha,\beta\right)\right) - \vard\left(\widehat{R}_{\frac{MC}{MC}}\right) = \frac{m}{n+m}\left( \vard\left(\widehat{R}_{\frac{CV}{CV}}\left(\alpha,\beta\right)\right)  - \vard\left(\widehat{R}_{\frac{MC}{MC}}\right)\right).
\end{align}

\newpage
\section{Additional results on the CV/CV estimator with linearly correlated control variates}
\label{sec:prop_B=D}
\subsection{Definition of alternative optimal coefficients}
The mild condition $|\corr(B,D)| < 1$ in Proposition \ref{prop_cv_cv}, which ensures that the control variates are not perfectly linearly correlated, can be easily satisfied in practice. It is sufficient to guarantee variance reduction using the optimal coefficients in Proposition \ref{prop_cv_cv}. While it is possible to use two linearly correlated (or even identical) control variates, in such cases the variance-minimizing coefficients $\Tilde{\alpha_o}$ and $\Tilde{\beta_o}$ defined in Proposition \ref{prop_B=D} no longer guarantee variance reduction, which is instead contingent upon the inequality \eqref{cond}.

\begin{proposition} \label{prop_B=D}
    Assume for all $(a,b)\in\mathbb{R}^*\times\mathbb{R}, B=aD+b$. The optimal coefficients for the CV/CV estimator with linearly correlated control variates $B$ and $D$ are
    \begin{align}
    \begin{cases}
    \Tilde{\alpha_o} = \frac{1}{a} \left( \frac{\cov(A,D) - R \cov(C,D)}{\var(D)} + \beta_o R \right); \\
    \Tilde{\beta_o} \in \mathbb{R},
    \end{cases}
    \end{align}
    where $R$ is the ratio to estimate.\\
    As a consequence, the estimator $ \widehat{R}_{\frac{CV}{CV}}\left(\Tilde{\alpha_o},\Tilde{\beta_o}\right)$ achieves variance reduction if the following condition holds:\begin{align}\cov(A-RC, D) \in \left(-\sqrt{\frac{\var(D)}{2}}, 0 \right) \cup \left(\sqrt{\frac{\var(D)}{2}}, +\infty \right). \label{cond}\end{align}
\end{proposition}
Details of the derivation are given in the following section. The previous condition can be checked to determine whether the CV/CV estimator achieves variance reduction.

\subsection{Proof}
\begin{proof}[Proof of Proposition \ref{prop_B=D}]
The setting where $\corr(B,D)=1$ is studied in this section. The two control variates $B$ and $D$ are linearly correlated if there is a linear relationship between them, i.e. $B=aD+b$ where $(a,b)\in\mathbb{R}^*\times\mathbb{R}$.

Some terms in the variance expression involving $B$ can be written in terms of $D$:
\begin{align*}
    \E[B] &= a\E[D] + b;\\
    \var(B) &= a^2\var(D); \\
    \cov(X,B) &= a\cov(X,D) \quad \forall X\in\{A,C\}.
\end{align*}

The CV/CV estimator's approximated variance can be expressed as:
\begin{align*}
\vard\left(\widehat{R}_{\frac{CV}{CV}}\left(\alpha,\beta\right)\right) =& \frac{1}{n\E[C]^2}
\left.
\begin{aligned}[t]
\Big( &\var(A) + a^2\alpha^2\var(D) - 2a\alpha\cov(A,D) + R^2\var(C)\\
&+ \beta^2R^2\var(D) - 2\beta R^2\cov(C,D) - 2R\cov(A,C) \\
&+ 2a\alpha R \cov(C,D) + 2\beta R \cov(A,D) - 2a\alpha\beta R\var(D) \Big) 
\end{aligned}
\right. \\
&+ \mathcal{O}_{\infty}(n^{-3/2}).
\end{align*}

This approximated variance is minimized by differentiating its expression with respect to $\alpha$ and $\beta$, and solving for the optimal value:
\begin{align*}
    \nabla \vard\left(\widehat{R}_{\frac{CV}{CV}}\left(\alpha,\beta\right)\right) &= \begin{bmatrix}
    \frac{\partial }{\partial \alpha} \vard\left(\widehat{R}_{\frac{CV}{CV}}\left(\alpha,\beta\right)\right)\\
    \frac{\partial }{\partial \beta} \vard\left(\widehat{R}_{\frac{CV}{CV}}\left(\alpha,\beta\right)\right)
    \end{bmatrix} \\
    &= \frac{2}{n\E[C]^2} \begin{bmatrix}
     a^2 \alpha \var(D) - a \cov(A,D) + a R \cov(C,D) - a \beta R \var(D) \\
     \beta R^2 \var(D) - R^2 \cov(C,D) + R \cov(A,D) - a \alpha R \var(D) 
    \end{bmatrix}.
\end{align*}

$\nabla \vard\left(\widehat{R}_{\frac{CV}{CV}}\left(\Tilde{\alpha_o},\Tilde{\beta_o}\right)\right) = \begin{bmatrix} 0 \\ 0 \end{bmatrix}$ is solved to find the expression of the minimizers $\Tilde{\alpha_o}$ and $\Tilde{\beta_o}$.
\begin{align*}
\nabla \vard\left(\widehat{R}_{\frac{CV}{CV}}\left(\Tilde{\alpha_o},\Tilde{\beta_o}\right)\right) = \begin{bmatrix} 0 \\ 0 \end{bmatrix} &\Leftrightarrow \begin{cases}
     \var(D) \left(a \Tilde{\alpha_o} - \Tilde{\beta_o} R\right) = \cov(A,D) - R \cov(C,D) \\
     \var(D) \left(\Tilde{\beta_o} R - a \Tilde{\alpha_o} \right) = R \cov(C,D) - \cov(A,D) 
\end{cases} \\
&\Leftrightarrow \Tilde{\alpha_o} = \frac{1}{a} \left( \frac{\cov(A,D) - R \cov(C,D)}{\var(D)} + \Tilde{\beta_o} R \right).
\end{align*}

If the control variates $B$ and $D$ are linearly correlated, there is an infinity of minimizers defined as \begin{align}
\begin{cases}
\Tilde{\alpha_o} = \frac{1}{a} \left( \frac{\cov(A,D) - R \cov(C,D)}{\var(D)} + \Tilde{\beta_o} R \right); \\
\Tilde{\beta_o} \in \mathbb{R}.
\end{cases}
\end{align}

Note $\eta := \frac{\cov(A,D) - R \cov(C,D)}{\var(D)}$. By replacing $\alpha$ and $\beta$ with their optimal values $\Tilde{\alpha_o}$ and $\Tilde{\beta_o}$, the approximated variance expression is as follows:
\begin{align*}
&\vard\left(\widehat{R}_{\frac{CV}{CV}}\left(\Tilde{\alpha_o},\Tilde{\beta_o}\right)\right) \\
=& \frac{1}{n\E[C]^2} 
\left.
\begin{aligned}[t]
\Big( &\var(A) + \left( \eta + \Tilde{\beta_o} R \right)^2\var(D) - 2\left( \eta + \Tilde{\beta_o} R \right)\cov(A,D) + R^2\var(C)\\
& + \Tilde{\beta_o}^2 R^2 \var(D) - 2\Tilde{\beta_o} R^2\cov(C,D) - 2R\cov(A,C) + 2\left( \eta + \Tilde{\beta_o} R \right) R \cov(C,D) \\
& + 2\Tilde{\beta_o} R \cov(A,D) - 2\left( \eta + \Tilde{\beta_o} R \right)\Tilde{\beta_o} R\var(D) \Big) + \mathcal{O}_{\infty}(n^{-3/2})
\end{aligned}
\right. \\
=& \frac{1}{n\E[C]^2}
\left.
\begin{aligned}[t]
\Big( &\var(A) + \eta^2 \var(D) - 2 \eta \cov(A,D) + R^2\var(C) - 2R\cov(A,C) \\
&+ 2 \eta R \cov(C,D) + 2 \eta \Tilde{\beta_o} R \var(D) + \Tilde{\beta_o}^2 R^2 \var(D) - 2 \Tilde{\beta_o} R \cov(A,D)\\
&+ \Tilde{\beta_o}^2 R^2 \var(D) - 2\Tilde{\beta_o} R^2\cov(C,D) + 2 \Tilde{\beta_o} R^2 \cov(C,D) + 2\Tilde{\beta_o} R \cov(A,D)\\
& - 2 \eta \Tilde{\beta_o} R\var(D) - 2 \Tilde{\beta_o}^2 R^2 \var(D) \Big) + \mathcal{O}_{\infty}(n^{-3/2})
\end{aligned}
\right. \\
=& \frac{1}{n\E[C]^2} \begin{aligned}[t] \Big(&\var(A) + \eta \var(D) - 2 \eta \cov(A,D) + R^2\var(C) - 2R\cov(A,C)\\
&+ 2 \eta R \cov(C,D) \Big) + \mathcal{O}_{\infty}(n^{-3/2}).
\end{aligned}
\end{align*}
The minimized variance is the same for any value of $\Tilde{\beta_o}$. Note that it is independent of $a$ and $b$, which means the reduced variance is the same for any control variates $B$ and $D=aB+d$.

The setting with two equal control variates $B=D$ is included in this study with $a=1$ and $b=0$.

The approximated variance reduction with linearly correlated control variates and optimal coefficients takes the following form:
\begin{align*}
&\vard\left(\widehat{R}_{\frac{CV}{CV}}\left(\Tilde{\alpha_o},\Tilde{\beta_o}\right)\right) - \vard\left(\widehat{R}_{\frac{MC}{MC}}\right) \\
=& \frac{1}{n\E[C]^2} \left(\eta \var(D) - 2 \eta \cov(A,D) + 2 \eta R \cov(C,D) \right) + \mathcal{O}_{\infty}(n^{-3/2})\\
=& \frac{1}{n\E[C]^2 \var(D)} \begin{aligned}[t]\Big( &\cov(A,D)\var(D) - R \cov(C,D)\var(D) - 2 \cov(A,D)^2 \\ &+ 4 R \cov(C,D) \cov(A,D) - 2 R^2 \cov(C,D)^2 \Big) + \mathcal{O}_{\infty}(n^{-3/2}) \end {aligned}\\
=& \frac{1}{n\E[C]^2} \left( \cov(A,D) - R \cov(C,D) - 2 \frac{\left(\cov(A,D) - R \cov(C,D)\right)^2}{\var(D)} \right) + \mathcal{O}_{\infty}(n^{-3/2})\\
=& \frac{1}{n\E[C]^2} \tau \left( 1 - 2 \frac{\tau^2}{\var(D)} \right) + \mathcal{O}_{\infty}(n^{-3/2})\\
&\text{where } \tau=\cov(A,D) - R \cov(C,D) = \cov(A-RC, D).
\end{align*}
The variance is reduced if $\tau \left( 1 - 2 \frac{\tau^2}{\var(D)} \right)<0$. That is
\begin{align*}
\begin{cases}
    \tau<0 \text{ and } \left( 1 - 2 \frac{\tau^2}{\var(D)} \right)>0 ;\\
    \tau>0 \text{ and } \left( 1 - 2 \frac{\tau^2}{\var(D)} \right)<0 .
\end{cases}
\quad \Leftrightarrow \quad \begin{cases}
    -\sqrt{\frac{\var(D)}{2}} < \cov(A-RC, D) < 0 ;\\
    \sqrt{\frac{\var(D)}{2}} < \cov(A-RC, D) .
\end{cases}
\end{align*}
Unlike the more general case where the control variates are non linearly correlated, when $B=aD+b$, the variance reduction is not guaranteed unless the condition $\cov(A-RC, D) \in \left(-\sqrt{\frac{\var(D)}{2}}, 0 \right) \cup \left(\sqrt{\frac{\var(D)}{2}}, +\infty \right) $ is fulfilled. \\
\end{proof}

\newpage
\section{Additional results on the CV/MC estimator} \label{sec:def_cv_mc}
\subsection{Definition of the CV/MC estimator}
\subsubsection{Exact control variates}
An alternative strategy to derive optimal coefficients for the CV/CV estimator is to focus on the numerator only. When the denominator variance is small compared to the numerator's variance, specifically when $R^2 \var(C) \ll \var(A)$, reducing it may be unnecessary.  In such cases, the denominator coefficient $\beta$ can be set to zero, and only the numerator coefficient $\alpha$ is optimized. We refer to the resulting estimator as the CV/MC estimator, known in the literature as the top-controlled estimator \cite{bauer1987control}, as the control variates method is applied only to the numerator, while the denominator is estimated by Monte Carlo.
\begin{definition}[CV/MC estimator] \label{def:CV/MC}
Let $(A_i, B_i, C_i)_{i=1 \dots n}$ be i.i.d. samples from the joint distribution of the random variables $A, B, C \in L^2$, where $B$ is a control variate. Assume $\E[B]$ is known. The CV/MC estimator of $R$ is defined as: \begin{align}\widehat{R}_{\frac{CV}{MC}}\left(\alpha_o'\right) = \frac{\widehat{A}_{CV}(\alpha_o')}{\overline{C_n}} = \frac{\overline{A_n} + \alpha_o' (\E[B] - \overline{B_n})}{\overline{C_n}},\end{align} where the coefficient $\alpha_o'$ is set to minimize the variance of the ratio as \begin{align}\alpha_o' := \underset{\alpha \in \mathbb{R}}{\text{argmin}} \; \vard\left(\frac{\overline{A_n} + \alpha (\E[B] - \overline{B_n})}{\overline{C_n}}\right) = \frac{\cov(A,B)- R\cov(B,C)}{\var(B)}. \label{eq:CV_MC}\end{align}
\end{definition}
Appendix \ref{sec:cv_mc} details the derivation of the optimal coefficient $\alpha_o'$.

\begin{proposition}
    With the optimal coefficient $\alpha_o'$ given in \eqref{eq:CV_MC}, the approximated variance reduction obtained with the CV/MC estimator is: \begin{align}
    \vard\left(\widehat{R}_{\frac{CV}{MC}}\left(\alpha_o'\right)\right) - \vard\left(\widehat{R}_{\frac{MC}{MC}}\right) &= - \frac{1}{n\E[C]^2} \frac{\left(\cov(A,B)-R\cov(B,C)\right)^2}{\var(B)} + \mathcal{O}_{\infty}(n^{-3/2}) \leq 0.
    \end{align}
    The variance reduction is guaranteed.
\end{proposition}
As the control variates mean estimator, the CV/MC ratio of means estimator presented here guarantees a variance no greater than that of the MC/MC estimator. 

The CV/MC estimator requires only the numerator coefficient $\alpha_o'$, which is simpler to estimate than the full pair $(\alpha_o, \beta_o)$ of the CV/CV estimator. Introducing the denominator coefficient $\beta$ adds an additional degree of freedom, allowing for further variance reduction at the cost of increased estimation complexity. Theoretically, the CV/CV estimator achieves a lower variance than the CV/MC estimator, as it minimizes the variance with respect to both $\alpha$ and $\beta$, whereas the CV/MC estimator restricts $\beta$ to zero. In practice, the CV/MC estimator may only perform comparably when the optimal $\beta_o$ is close to zero.

\subsubsection{Approximate control variates}
The CV/MC estimator previously introduced can be easily extended to the approximate control variates method, by replacing the control variates mean $\E[B]$ by Monte Carlo estimation on $n+m$ points:
\begin{align}
    \widehat{R}_{\frac{ACV}{MC}}\left(\alpha\right) = \frac{\overline{A_n} + \alpha (\overline{B_{n+m}} - \overline{B_n})}{\overline{C_n}}.
\end{align}
A similar approximated variance reduction formula holds for the ACV/MC estimator, which scales the exact CV/MC case by a factor of $\frac{m}{n+m}$: \begin{align}
    \vard\left(\widehat{R}_{\frac{ACV}{MC}}\left(\alpha\right)\right) - \vard\left(\widehat{R}_{\frac{MC}{MC}}\right) = \frac{m}{n+m} \left( \vard\left(\widehat{R}_{\frac{CV}{MC}}\left(\alpha\right)\right) - \vard\left(\widehat{R}_{\frac{MC}{MC}}\right) \right).
\end{align}

\subsection{Simulation study}
The reported simulation values are averaged over 10,000 independent repetitions, comparing different coefficients. For the CV/MC estimator, the coefficient considered is  $\alpha_o'$ the optimal coefficient, given in \eqref{eq:CV_MC}, that minimizes the variance of the full ratio.\\
For the CV/CV estimator, three sets of coefficients are considered:
\begin{itemize}
    \item the "classical" coefficients $(\alpha_c, \beta_c)$, given in \eqref{eq:classic_alpha} and \eqref{eq:classic_beta}, derived by optimizing the numerator and denominator independently;
    \item the  coefficients proposed by \cite{gordon_efficient_1982} $(\alpha_g, \beta_g)$, given in \eqref{eq:gordon_alpha} and \eqref{eq:gordon_beta}, where $\alpha_g$ minimizes the numerator variance and $\beta_g$ minimizes the approximated variance of the ratio;
    \item the new optimal coefficients $(\alpha_o, \beta_o)$, given in \eqref{eq:alpha_opt} and \eqref{eq:beta_opt}, that jointly minimize the approximated variance of the ratio.
\end{itemize}

\begin{figure}[!h]
    \centering
    \includegraphics[height=0.4\textheight]{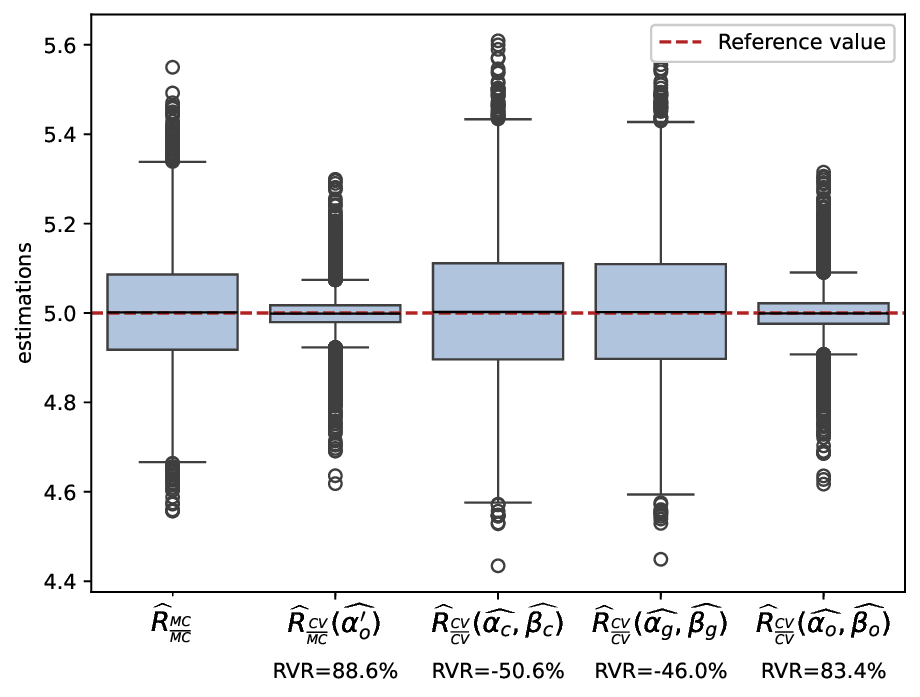}
    \vspace{-10pt}
    \caption{Boxplots of ratio of means estimations, comparing the MC/MC, CV/MC and CV/CV estimators - Variance increase with non-optimal coefficients.}
    \label{fig:boxplot_cv_mc}
\end{figure}
The presented results are observed with $n=10$ samples and the covariance matrix for the worst-case scenario for the coefficients of \cite{gordon_efficient_1982}, as defined in Section \ref{sec:worst_gordon}.

Figure \ref{fig:boxplot_cv_mc} demonstrates that only the optimal coefficients consistently guarantee variance reduction. In this specific scenario, the CV/MC estimator outperforms the CV/CV approach by 5\% of RVR; this is attributed to the low correlation of the control variate $D$ with the other variables. Furthermore, the higher complexity of the CV/CV optimal coefficients leads to estimation instability at small sample sizes. However, as the sample size increases to $n=100$, the performance gap vanishes, suggesting that the relative advantage of the simpler CV/MC estimator is primarily a finite-sample effect observed in low-correlation regimes.

The measured RAE of 9.4\% suggests the approximation of the variance remains acceptable.

\subsection{Proofs} \label{sec:cv_mc}
\subsubsection{Approximated variance of the CV/MC estimator}
An expression of the CV/MC estimator's variance is derived, so that the optimal coefficient minimizing it can be found.

To derive the CV/MC estimator's variance expression, a few preliminary results are needed: 
\begin{lemma} \label{lem:cv_mc}
    \begin{enumerate}
        \item $\E[\widehat{A}_{CV}(\alpha)] = \E[A]$;
        \item $\var(\widehat{A}_{CV}(\alpha)) = \frac{1}{n}\var(A) + \frac{\alpha^2}{n} \var(B) -2\frac{\alpha}{n} \cov(A, B)$;
        \item $\cov(\widehat{A}_{CV}(\alpha), \overline{C_n}) = \frac{1}{n}\cov(A, C) - \frac{\alpha}{n} \cov(B, C)$.
    \end{enumerate}
\end{lemma}
\begin{proof}
    Point 1 can be proven as \begin{align*}\E[\widehat{A}_{CV}(\alpha)] = \E[\overline{A_n} + \alpha (\E[B] - \overline{B_n})] = \E[\overline{A_n}] + \alpha (\E[B] - \E[\overline{B_n}]) = \E[A]. \end{align*}
    Point 2 can be proven as \begin{align*}
    \var(\widehat{A}_{CV}(\alpha)) &= \var(\overline{A_n} + \alpha(\E[B] - \overline{B_n})) \\
    &= \var(\overline{A_n}) + \alpha^2 \var(\E[B] - \overline{B_n}) +2\alpha \cov(\overline{A_n}, \E[B] - \overline{B_n}) \\
    &= \var(\overline{A_n}) + \alpha^2 \var(\overline{B_n}) -2\alpha \cov(\overline{A_n}, \overline{B_n}) \\
    &= \frac{1}{n}\var(A) + \frac{\alpha^2}{n} \var(B) -2\frac{\alpha}{n} \cov(A, B).
    \end{align*}
    Point 3 can be proven as \begin{align*}
    \cov(\widehat{A}_{CV}(\alpha), \overline{C_n}) &= \cov(\overline{A_n} + \alpha(\E[B] - \overline{B_n}), \overline{C_n}) \\
    &= \cov(\overline{A_n}, \overline{C_n}) + \alpha \cov(\E[B] - \overline{B_n}, \overline{C_n}) \\
    &= \cov(\overline{A_n}, \overline{C_n}) - \alpha \cov(\overline{B_n}, \overline{C_n}) \\
    &= \frac{1}{n}\cov(A, C) - \frac{\alpha}{n} \cov(B, C).
    \end{align*}
\end{proof}

\begin{proposition}[Variance of the CV/MC estimator]
Let $(A_i, B_i, C_i)_{i=1 \dots n}$ be i.i.d. samples from the joint distribution of the random variables $A, B, C \in L^2$, where $B$ is a control variate. The variance of the CV/MC estimator can be approximated as
\begin{equation}
\boxed{
\begin{aligned}[t]
\vard\left(\widehat{R}_{\frac{CV}{MC}}\left(\alpha\right)\right) = \frac{1}{n\E[C]^2} \bigg( &\var(A) + \alpha^2\var(B) - 2\alpha\cov(A, B) + R^2\var(C) \\ 
&- 2R\cov(A, C) + 2\alpha R\cov(B, C) \bigg) + \mathcal{O}_{\infty}(n^{-3/2}).
\end{aligned}
} \label{eq:var_cv_mc}
\end{equation}
\end{proposition}

\begin{proof}
    Using the expressions in Lemma \ref{lem:cv_mc} and Proposition \ref{prop:var_mc_mc}, the variance of the CV/MC estimator can be approximated as:
    
    \begin{align*}
        \vard\left(\widehat{R}_{\frac{CV}{MC}}\left(\alpha\right)\right)
        =&  \frac{1}{\E[\overline{C_n}]^2}\var(\widehat{A}_{CV}(\alpha)) + \frac{\E[\widehat{A}_{CV}(\alpha)]^2}{\E[\overline{C_n}]^4} \var(\overline{C_n}) \\
        &-2 \frac{\E[\widehat{A}_{CV}(\alpha)]}{\E[\overline{C_n}]^3} \cov(\widehat{A}_{CV}(\alpha), \overline{C_n}) + \mathcal{O}_{\infty}(n^{-3/2})\\ 
        =&  \frac{1}{\E[C]^2}\var(\widehat{A}_{CV}(\alpha)) + \frac{\E[A]^2}{\E[C]^4} \var(\overline{C_n}) \\&-2 \frac{\E[A]}{\E[C]^3} \cov(\widehat{A}_{CV}(\alpha), \overline{C_n}) + \mathcal{O}_{\infty}(n^{-3/2})\\ 
        =&  \frac{1}{\E[C]^2}\var(\overline{A_n}) + \frac{\alpha^2}{\E[C]^2}\var(\overline{B_n}) - \frac{2\alpha}{\E[C]^2}\cov(\overline{A_n}, \overline{B_n}) + \frac{\E[A]^2}{\E[C]^4} \var(\overline{C_n}) \\ 
        & -2 \frac{\E[A]}{\E[C]^3} \cov(\overline{A_n}, \overline{C_n}) + 2\alpha \frac{\E[A]}{\E[C]^3} \cov(\overline{B_n}, \overline{C_n}) + \mathcal{O}_{\infty}(n^{-3/2})\\ 
        =&  \frac{1}{\E[C]^2}\frac{1}{n}\var(A) + \frac{\alpha^2}{\E[C]^2}\frac{1}{n}\var(B) - \frac{2\alpha}{\E[C]^2}\frac{1}{n}\cov(A, B) + \frac{R^2}{\E[C]^2} \frac{1}{n} \var(C)\\ 
        & -2 \frac{R}{\E[C]^2} \frac{1}{n} \cov(A, C) + 2\alpha \frac{R}{\E[C]^2} \frac{1}{n} \cov(B, C) + \mathcal{O}_{\infty}(n^{-3/2}).
    \end{align*}
\end{proof}

\vspace{10pt}
\subsubsection{Optimal coefficient of the CV/MC estimator}
\begin{proof}[Proof of Definition \ref{def:CV/MC}]
    The goal is to find the optimal coefficient $\alpha_o'$ minimizing the variance of the CV/MC estimator.
    \begin{align*}
        \alpha_o' := \underset{\alpha \in \mathbb{R}}{\text{argmin}} \; \vard\left(\frac{\overline{A_n} + \alpha (\E[B] - \overline{B_n})}{\overline{C_n}}\right). \\
    \end{align*}
    The variance in \eqref{eq:var_cv_mc} is minimized by differentiating its expression with respect to $\alpha$ and solving for the optimal value:
    \begin{align*}
        \frac{d}{d \alpha}\vard\left(\widehat{R}_{\frac{CV}{MC}}\left(\alpha\right)\right) =  \frac{2 \alpha}{\E[C]^2}\var(\overline{B_n}) - \frac{2}{\E[C]^2}\cov(\overline{A_n}, \overline{B_n}) + 2 \frac{\E[A]}{\E[C]^3} \cov(\overline{B_n}, \overline{C_n}) + \mathcal{O}_{\infty}(n^{-3/2}) .
    \end{align*}
    The function $\alpha \mapsto \vard\left(\widehat{R}_{\frac{CV}{MC}}\left(\alpha\right)\right)$ is strictly convex as $\frac{d^2}{d \alpha^2} \vard\left(\widehat{R}_{\frac{CV}{MC}}\left(\alpha\right)\right) = \frac{2}{\E[C]^2}\var(\overline{B_n}) > 0$, so the value of $\alpha$ for which $\frac{d}{d \alpha}\vard\left(\widehat{R}_{\frac{CV}{MC}}\left(\alpha\right)\right) = 0$ is the unique minimizer of the approximated variance.

    It gives:
    \begin{align*}
        \frac{d}{d \alpha}\vard\left(\widehat{R}_{\frac{CV}{MC}}\left(\alpha_o'\right)\right) = 0 \Leftrightarrow \alpha_o' &= \frac{\E[C]^2}{2 \var(\overline{B_n})} \left(\frac{2}{\E[C]^2}\cov(\overline{A_n}, \overline{B_n}) - 2 \frac{\E[A]}{\E[C]^3} \cov(\overline{B_n}, \overline{C_n})\right)\\ &= \frac{\cov(\overline{A_n}, \overline{B_n})}{\var(\overline{B_n})} - \frac{\E[A] \cov(\overline{B_n}, \overline{C_n})}{\E[C] \var(\overline{B_n})}.
    \end{align*}
    \begin{align*}
    \boxed{
    \alpha_o' = \frac{\cov(A, B)}{\var(B)} - R\frac{\cov(B, C)}{\var(B)}
    }
    \end{align*}
\end{proof}

\vspace{10pt}
\subsubsection{Comparing the approximated variances of the MC/MC and the CV/MC estimators with different coefficients}
The difference between approximated variances can be expressed as: \begin{align*}
&\vard\left(\widehat{R}_{\frac{CV}{MC}}\left(\alpha\right)\right) - \vard\left(\widehat{R}_{\frac{MC}{MC}}\right) \\
=& \frac{1}{n\E[C]^2} \left( \var(A) + \alpha^2\var(B) - 2\alpha\cov(A, B) + R^2\var(C) - 2R\cov(A, C) + 2\alpha R\cov(B, C) \right)\\
&- \frac{1}{n\E[C]^2} \left( \var(A) + R^2 \var(C) -2 RCov(A, C) \right) + \mathcal{O}_{\infty}(n^{-3/2})\\
=& \frac{1}{n \E[C]^2} \left( \alpha^2\var(B) - 2\alpha\cov(A, B) -2R\cov(A, C) + 2\alpha R \cov(B, C) \right) + \mathcal{O}_{\infty}(n^{-3/2}).
\end{align*}

$\bullet$ With $\alpha_c=\frac{\cov(A,B)}{\var(B)}$ the classical coefficient, the approximated variance difference is \begin{align*}
\vard\left(\widehat{R}_{\frac{CV}{MC}}\left(\alpha_c\right)\right) - \vard\left(\widehat{R}_{\frac{MC}{MC}}\right) = \frac{1}{n\E[C]^2} \frac{\cov(A,B)}{\var(B)} \left( 2R \cov(B,C) - \cov(A,B) \right) + \mathcal{O}_{\infty}(n^{-3/2}).
\end{align*}
So there is variance reduction if and only if \\$\vard\left(\widehat{R}_{\frac{CV}{MC}}\left(\alpha_c\right)\right) - \vard\left(\widehat{R}_{\frac{MC}{MC}}\right) <0 \Leftrightarrow \cov(A,B) \left( 2R \cov(B,C) - \cov(A,B) \right) <0$. \\

$\bullet$ With $\alpha_o' = \frac{\cov(A,B)}{\var(B)} - R\frac{\cov(B,C)}{\var(B)}$ the optimized coefficient, the approximated variance difference is \begin{align*}
\vard\left(\widehat{R}_{\frac{CV}{MC}}\left(\alpha_o'\right)\right) - \vard\left(\widehat{R}_{\frac{MC}{MC}}\right) = - \frac{1}{n\E[C]^2} \frac{\left(\cov(A,B)-R\cov(B,C)\right)^2}{\var(B)} + \mathcal{O}_{\infty}(n^{-3/2}) \leq 0. 
\end{align*}
The variance reduction is guaranteed.

\end{document}